# ASYMPTOTIC OPERATING CHARACTERISTICS OF AN OPTIMAL CHANGE POINT DETECTION IN HIDDEN MARKOV MODELS

By Cheng-Der Fuh[1]

*Academia Sinica*

Let $\xi_0, \xi_1, \ldots, \xi_{\omega-1}$ be observations from the hidden Markov model with probability distribution $P^{\theta_0}$, and let $\xi_\omega, \xi_{\omega+1}, \ldots$ be observations from the hidden Markov model with probability distribution $P^{\theta_1}$. The parameters $\theta_0$ and $\theta_1$ are given, while the change point $\omega$ is unknown. The problem is to raise an alarm as soon as possible after the distribution changes from $P^{\theta_0}$ to $P^{\theta_1}$, but to avoid false alarms. Specifically, we seek a stopping rule $N$ which allows us to observe the $\xi's$ sequentially, such that $E_\infty N$ is large, and subject to this constraint, $\sup_k E_k(N-k|N \geq k)$ is as small as possible. Here $E_k$ denotes expectation under the change point $k$, and $E_\infty$ denotes expectation under the hypothesis of no change whatever.

In this paper we investigate the performance of the Shiryayev–Roberts–Pollak (SRP) rule for change point detection in the dynamic system of hidden Markov models. By making use of Markov chain representation for the likelihood function, the structure of asymptotically minimax policy and of the Bayes rule, and sequential hypothesis testing theory for Markov random walks, we show that the SRP procedure is asymptotically minimax in the sense of Pollak [*Ann. Statist.* **13** (1985) 206–227]. Next, we present a second-order asymptotic approximation for the expected stopping time of such a stopping scheme when $\omega = 1$. Motivated by the sequential analysis in hidden Markov models, a nonlinear renewal theory for Markov random walks is also given.

**1. Introduction.** The problem of quick detection, with low false-alarm rate, of abrupt changes in stochastic dynamic systems arises in a variety of

Received October 2003; revised March 2004.
[1]Supported in part by NSC 92-2118-M-001-032.

applications, including industrial quality control, segmentation of signals, financial engineering, biomedical signal processing, edge detection in images, and the diagnosis of faults in the elements of computer communication networks. A comprehensive summary in this area was given by Basseville and Nikiforov (1993) and Lai (1995, 2001). A typical such problem in segmentation of signals is that of using an automatic segmentation of the signal as the first processing step, and a segmentation algorithm splits the signal into homogeneous segments, the lengths of which are adapted to the local characteristics of the analyzed signal. The main desired properties of a segmentation algorithm are few false alarms and missed detections, and low detection delay. In the standard formulation of the change point detection problem, there is a sequence of observations whose distribution changes at some unknown time $\omega$, and the goal is to detect this change as soon as possible under false alarm constraints. The reader is referred to Braun and Müller (1998) for a nice discussion of hidden Markov models for DNA data and change point detection analysis.

When the observations $\xi_n$ are independent with a common density function $f^{\theta_0}$ for $n < \omega$ and with another common density function $f^{\theta_1}$ for $n \geq \omega$, a minimax formulation has been proposed by Lorden (1971), in which he showed that subject to the "average run length" (ARL) constraint, Page's CUSUM procedure asymptotically minimizes the "worst case" detection delay. Instead of studying the optimal detection problem via sequential testing theory, Moustakides (1986) formulated the worst case detection delay problem subject to an ARL constraint as an optimal solution to the optimal stopping problem. Ritov (1990) later gave a simpler proof. For change point detection in complex dynamic systems beyond the i.i.d. setting, Bansal and Papantoni-Kazakos (1986) extended Lorden's asymptotic theory to the case where $\xi_j$ are stationary ergodic sequences, under the condition that $\{\xi_j, j < \omega\}$ (before the change point) and $\{\xi_j, j \geq \omega\}$ (after the change point) are independent, and proved the asymptotic optimality of the CUSUM algorithm. Further extensions to general stochastic sequences $\xi_n$ were obtained by Lai (1995, 1998). Moreover, using a change of measure argument, Lai (1998) also established the asymptotic optimality of the CUSUM rule under several alternative performance criteria. In the dynamic system of hidden Markov models, Fuh (2003) proved that the CUSUM scheme is asymptotically optimal in the sense of Lorden (1971). His method related the CUSUM procedure to certain one-sided sequential probability ratio tests in hidden Markov models, for which they had been shown, in Section 4 of Fuh (2003), to be asymptotically optimal for testing simple hypotheses.

In the simple system of independent observations before and after the change, a Bayesian formulation has been proposed by Shiryayev (1963, 1978), in which the change point is assumed to have a geometric prior distribution, and the goal is to minimize the expected delay subject to an upper bound on the false alarm probability. He used optimal stopping theory to show that the Bayes rule triggers an alarm as soon as the posterior probability that a change has occured exceeds some fixed level. Roberts (1966) considered the non-Bayesian setting, and studied by simulation the average run length of this rule, and found it to be very good. Pollak and Siegmund (1975) extended Shiryayev's work in a non-Bayesian setting. And Pollak (1985) showed that the (modified) Shiryayev–Roberts rule is asymptotically minimax under the formulation of Pollak and Siegmund (1975). Later Yakir (1997) proved that the procedure is strictly optimal for a slight reformulation of the problem. Finally, we mention that Yakir (1994) studied Bayesian optimal detection for a finite state Markov chain.

As noted by Basseville and Nikiforov (1993) in their monograph, there is a great deal of literature on detection algorithms in complex systems but relatively little on the statistical properties and optimality theory of detection procedures beyond very simple models. The primary goal of this paper is to investigate theoretical aspects of the Shiryayev–Roberts–Pollak (SRP) change point detection rule in hidden Markov models. We show that the SRP procedure is asymptotically minimax in the sense of Pollak (1985). Next, we present a second-order asymptotic approximation for the expected stopping time of such a stopping scheme when $\omega = 1$. Motivated by the sequential analysis in hidden Markov models, a nonlinear renewal theory for Markov random walks is also given.

This paper is organized as follows. In Section 2 we define the hidden Markov model and formulate the sequential change point detection problem. Then we provide a Markov chain representation of the likelihood ratio. A nonlinear Markov renewal theory is given in Section 3. In Section 4 we show that the SRP rule is asymptotically minimax under mild conditions. In Section 5 we study the asymptotic operating characteristics of the detection procedure, and derive a second-order asymptotic approximation for the expected stopping scheme when $\omega = 1$. All proofs are given in Sections 6, 7 and 8.

**2. Problem formulation.** A hidden Markov model is defined as a parameterized Markov chain in a Markovian random environment [Fuh (2003)], with the underlying environmental Markov chain viewed as missing data. That is, for each $\theta \in \Theta \subset R^q$, the unknown parameter, we consider $\mathbf{X} = \{X_n, n \geq 0\}$ as an ergodic (positive recurrent, irreducible and aperiodic) Markov chain on a finite state space $D = \{1, 2, \ldots, d\}$, with transition probability matrix $P(\theta) = [p_{xy}(\theta)]_{x,y=1,\ldots,d}$ and stationary distribution $\pi(\theta) = (\pi_x(\theta))_{x=1,\ldots,d}$. Suppose that an additive component $\xi_n$, taking values in $R$, is adjoined to the chain such that $\{(X_n, \xi_n), n \geq 0\}$ is a Markov chain on $D \times R$, satisfying $P^{(\theta)}\{X_1 \in A | X_0 = x, \xi_0 = s\} = P^{(\theta)}\{X_1 \in A | X_0 = x\}$ for $A \in \mathcal{B}(D)$. And conditioning on the full **X** sequence, $\xi_n$ is a Markov chain with probability

$$
\begin{aligned}
P^\theta\{\xi_{n+1} &\in B | X_0, X_1, \ldots; \xi_0, \xi_1, \ldots, \xi_n\} \\
&= P^\theta\{\xi_{n+1} \in B | X_{n+1}; \xi_n\} = P^\theta(X_{n+1} : \xi_n, B) \qquad \text{a.s.}
\end{aligned}
\tag{2.1}
$$

for each $n$ and $B \in \mathcal{B}(R)$, the Borel $\sigma$-algebra of $R$. Note that in (2.1) the conditional probability of $\xi_{n+1}$ depends on $X_{n+1}$ and $\xi_n$ only. Furthermore, we assume the existence of a transition probability density for the Markov chain $\{(X_n, \xi_n), n \geq 0\}$ with respect to a $\sigma$-finite measure $\mu$ on $R$ such that

$$
\begin{aligned}
P^{(\theta)}&\{X_1 \in A, \xi_1 \in B | X_0 = x, \xi_0 = s_0\} \\
&= P^{(\theta)}\{\xi_1 \in B | X_1 \in A, X_0 = x, \xi_0 = s_0\} P^{(\theta)}\{X_1 \in A | X_0 = x\} \\
&= \sum_{y \in A} \int_B p_{xy}(\theta) f(s; \varphi_y(\theta) | s_0) \, d\mu(s),
\end{aligned}
\tag{2.2}
$$

where $f(\xi_k; \varphi_{X_k}(\theta) | \xi_{k-1})$ is the transition probability density of $\xi_k$ given $\xi_{k-1}$ and $X_k$ with respect to $\mu$, $\theta \in \Theta$ is the unknown parameter, and $\varphi_y(\cdot)$ is a function defined on the parameter space $\Theta$ for each $y = 1, \ldots, d$. Here and in the sequel we assume the Markov chain $\{(X_n, \xi_n), n \geq 0\}$ has stationary probability $\Gamma$ with probability density $\pi_x(\theta) f(\cdot; \varphi_x(\theta))$ with respect to $\mu$. In this paper we assume that only one parameter is of interest and treat the other parameters as nuisance parameters. That is, for simplicity we consider $\theta \in \Theta \subseteq R$ as a one-dimensional unknown parameter. For convenience of notation, we write $\pi_x$ for $\pi_x(\theta)$ and $p_{xy}$ for $p_{xy}(\theta)$. We call a process $\{\xi_n, n \geq 0\}$ a hidden Markov model if there is a Markov chain $\{X_n, n \geq 0\}$ such that the process $\{(X_n, \xi_n), n \geq 0\}$ satisfies (2.1) and (2.2).

Let $\xi_0, \xi_1, \ldots, \xi_{\omega-1}$ be the observations from the hidden Markov model $\{\xi_n, n \geq 0\}$ with distribution $P^{\theta_0}$, and let $\xi_\omega, \xi_{\omega+1}, \ldots$ be the observations from the hidden Markov model $\{\xi_n, n \geq 0\}$ with distribution $P^{\theta_1}$. Both $\theta_0$ and $\theta_1$ are given, while the change point $\omega$ is unknown. We shall use $P_\omega$ to denote such a probability measure (with change time $\omega$) and use $P_\infty$ to denote the case $\omega = \infty$ (no change point). Denote $E_\omega$ as the corresponding expectation under $P_\omega$. The objectives are to raise an alarm as soon as possible after the change and to avoid false alarms. A detection scheme is a stopping time on the sequence of observations and aims to minimize the number of post change observations. Hence, the stopping time $N$ should satisfy $\{N \geq \omega\}$ but, at the same time, keep $N - \omega$ small. In this paper we use the functional studied by Pollak and Siegmund (1975) and Pollak (1985) to find a stopping time $N$ to minimize

$$
\sup_{1 \leq k < \infty} E_k(N - k | N \geq k)
\tag{2.3}
$$

subject to

(2.4) $$E_\infty N \geq \gamma,$$

for some specified (large) constant $\gamma$. A detection scheme is called asymptotically minimax if it minimizes (2.3), within an $o(1)$ order, among all stopping rules that satisfy $E_\infty N \geq \gamma$, where $o(1) \to 0$ as $\gamma \to \infty$.

To describe the SRP change point detection scheme, we need the following notation. Fix $\theta_0, \theta_1 \in \Theta$. Let $\xi_0, \xi_1, \ldots, \xi_n$ be the observations given from the hidden Markov model $\{\xi_n, n \geq 0\}$. Denote

(2.5)
$$\begin{aligned}
LR_n &:= \frac{p_n(\xi_0, \xi_1, \ldots, \xi_n; \theta_1)}{p_n(\xi_0, \xi_1, \ldots, \xi_n; \theta_0)} \\
&:= \sum_{x_0=1}^d \cdots \sum_{x_n=1}^d \pi_{x_0}(\theta_1) f(\xi_0; \varphi_{x_0}(\theta_1)) \\
&\quad \times \prod_{l=1}^n p_{x_{l-1}x_l}(\theta_1) f(\xi_l; \varphi_{x_l}(\theta_1)|\xi_{l-1}) \\
&\quad \times \left[ \sum_{x_0=1}^d \cdots \sum_{x_n=1}^d \pi_{x_0}(\theta_0) f(\xi_0; \varphi_{x_0}(\theta_0)) \right. \\
&\quad\quad \left. \times \prod_{l=1}^n p_{x_{l-1}x_l}(\theta_0) f(\xi_l; \varphi_{x_l}(\theta_0)|\xi_{l-1}) \right]^{-1}
\end{aligned}$$

as the likelihood ratio. For $0 \leq k \leq n$, let

(2.6)
$$\begin{aligned}
LR_n^k &:= \frac{p_n(\xi_k, \xi_{k+1}, \ldots, \xi_n; \theta_1)}{p_n(\xi_k, \xi_{k+1}, \ldots, \xi_n; \theta_0)} \\
&:= \frac{\sum_{x_k=1}^d \cdots \sum_{x_n=1}^d \prod_{l=k}^n p_{x_{l-1}x_l}(\theta_1) f(\xi_l; \varphi_{x_l}(\theta_1)|\xi_{l-1})}{\sum_{x_k=1}^d \cdots \sum_{x_n=1}^d \prod_{l=k}^n p_{x_{l-1}x_l}(\theta_0) f(\xi_l; \varphi_{x_l}(\theta_0)|\xi_{l-1})}.
\end{aligned}$$

Given an approximate threshold $B > 0$ and setting $b = \log B$, define the Shiryayev–Roberts scheme

(2.7) $$N_b := \inf\left\{n : \sum_{k=0}^n LR_n^k \geq B\right\} = \inf\left\{n : \log \sum_{k=0}^n LR_n^k \geq b\right\}.$$

A simple modification of (2.7) was given by Pollak (1985) by adding a randomization on the initial $LR_n^0$. This will be defined precisely in Section 4.

It is worth asking that while the SRP rule (2.5)–(2.7) is asymptotically minimax in the i.i.d. cases [Pollak (1985)], is it nontrivial whether this is still true for hidden Markov models? To give a definitive answer to this question, we need to study the likelihood ratio $LR_n$ that appeared in (2.5) since (2.6) can be analyzed in the same manner. Note that the nonadditive form of (2.5) makes it difficult to analyze. A key idea to get rid of this difficulty is to represent the likelihood ratio (2.5) as the ratio of $L_1$-norms of products of Markov random matrices. This device has been proposed by Fuh (2003) to study SPRT and CUSUM for hidden Markov models. Here, we carry out the same idea to have a representation of the likelihood ratio $LR_n$.

Given a column vector $u = (u_1, \ldots, u_d)^t \in R^d$, where $t$ denotes the transpose of the underlying vector in $R^d$, define the $L_1$-norm of $u$ as $\|u\| = \sum_{i=1}^d |u_i|$. The likelihood ratio (2.5) can be represented as

$$(2.8) \quad LR_n = \frac{p_n(\xi_0, \xi_1, \ldots, \xi_n; \theta_1)}{p_n(\xi_0, \xi_1, \ldots, \xi_n; \theta_0)} = \frac{\|M_n(\theta_1) \cdots M_1(\theta_1) M_0(\theta_1) \pi(\theta_1)\|}{\|M_n(\theta_0) \cdots M_1(\theta_0) M_0(\theta_0) \pi(\theta_0)\|}$$

where, for $\theta = \theta_0$ or $\theta_1$,

$$(2.9) \quad M_0 = M_0(\theta) = \begin{bmatrix} f(\xi_0; \varphi_1(\theta)) & 0 & \cdots & 0 \\ \vdots & \ddots & \vdots & \vdots \\ 0 & 0 & \cdots & f(\xi_0; \varphi_d(\theta)) \end{bmatrix},$$

$$(2.10) \quad M_k = M_k(\theta) = \begin{bmatrix} p_{11}(\theta) f(\xi_k; \varphi_1(\theta) | \xi_{k-1}) & \cdots & p_{d1}(\theta) f(\xi_k; \varphi_1(\theta) | \xi_{k-1}) \\ \vdots & \ddots & \vdots \\ p_{1d}(\theta) f(\xi_k; \varphi_d(\theta) | \xi_{k-1}) & \cdots & p_{dd}(\theta) f(\xi_k; \varphi_d(\theta) | \xi_{k-1}) \end{bmatrix}$$

for $k = 1, \ldots, n$, and

$$(2.11) \quad \pi(\theta) = (\pi_1(\theta), \ldots, \pi_d(\theta))^t.$$

Let $\{(X_n, \xi_n), n \geq 0\}$ be the Markov chain defined in (2.1) and (2.2). Denote $Y_n := (X_n, \xi_n)$ and $D' := D \times R$. Define $Gl(d, R)$ as the set of invertible $d \times d$ matrices with real entries. For given $k = 0, 1, \ldots, n$, and $\theta = \theta_0$ or $\theta_1$, let $M_k(\theta)$ be the random matrix from $D' \times D'$ to $Gl(d, R)$, as defined in (2.9) and (2.10). For convenience of notation, we still denote $\theta = (\theta_0, \theta_1)$ and let

$$(2.12) \quad \begin{aligned} \mathbb{T}_n(\theta) &= M_n(\theta) \cdots M_0(\theta) \\ &= (\mathbb{T}_n(\theta_0), \mathbb{T}_n(\theta_1)) = (M_n(\theta_0) \cdots M_0(\theta_0), M_n(\theta_1) \cdots M_0(\theta_1)). \end{aligned}$$

Then the system $\{(Y_n, \mathbb{T}_n(\theta)), n \geq 0\}$ is called a product of Markov random matrices on $D' \times Gl(d, R) \times Gl(d, R)$. Denote $\mathcal{P}_y^\theta$ as the probability distribution of $\{(Y_n, \mathbb{T}_n(\theta)), n \geq 0\}$ with $Y_0 = y$, and $\mathcal{E}_y^\theta$ as the expectation under $\mathcal{P}_y^\theta$.

Let $u \in R^d$ be a $d$-dimensional vector, $\bar{u} := u/\|u\|$ the normalization of $u$ ($\|u\| \neq 0$), and denote $P(R^d)$ as the projection space of $R^d$ which contains all elements $\bar{u}$. For given $\bar{u} \in P(R^d)$ and $M \in Gl(d, R)$, denote $M \cdot \bar{u} = \overline{Mu}$ and $\overline{\mathbb{T}_k(\theta)u} = (\overline{\mathbb{T}_k(\theta_0)u}, \overline{\mathbb{T}_k(\theta_1)u})$, for $k = 0, \ldots, n$. Let

$$(2.13) \quad W_0^\theta = (Y_0, \overline{\mathbb{T}_0(\theta)u}), \quad W_1^\theta = (Y_1, \overline{\mathbb{T}_1(\theta)u}), \ldots, W_n^\theta = (Y_n, \overline{\mathbb{T}_n(\theta)u}).$$

Then $\{W_n^\theta, n \geq 0\}$ is a Markov chain on the state space $D' \times P(R^d) \times P(R^d)$ with the transition kernel

$$(2.14) \qquad \mathbb{P}^\theta((y,\bar{u}), A \times B) := \mathcal{E}_y^\theta(I_{A \times B}(Y_1, \overline{M_1(\theta)u}))$$

for all $y \in D'$, $\bar{u} := (\bar{u}, \bar{u}) \in P(R^d) \times P(R^d)$, $A \in \mathcal{B}(\mathcal{D}')$, and $B \in \mathcal{B}(P(R^d) \times P(R^d))$, the Borel $\sigma$-algebra of $P(R^d) \times P(R^d)$. For simplicity we let $\mathbb{P}_{(y,\bar{u})}^\theta := \mathbb{P}^\theta(\cdot,\cdot)$ and denote $\mathbb{E}_{(y,\bar{u})}^\theta$ as the expectation under $\mathbb{P}_{(y,\bar{u})}^\theta$. Since the Markov chain $\{(X_n, \xi_n), n \geq 0\}$ has transition probability density and the random matrix $M_1(\theta)$ is driven by $\{(X_n, \xi_n), n \geq 0\}$, it implies that the induced transition probability $\mathbb{P}(\cdot,\cdot)$ has a density with respect to $\mu$. Denote it as $\mathbb{P}$ for simplicity. Under Condition C given below, the Markov chain $\{W_n^\theta, n \geq 0\}$ has an invariant probability measure $m^\theta$ on $D' \times P(R^d) \times P(R^d)$; see Fuh (2003).

Now, for $y_0, y_1 \in D'$, $\bar{u} = \overline{u(\theta)} = (\overline{u(\theta_0)}, \overline{u(\theta_1)}) \in P(R^d) \times P(R^d)$ and $M = M(y_0, y_1) = M(\theta) = (M(\theta_0), M(\theta_1)) \in Gl(d, R) \times Gl(d, R)$, let $\sigma : (D' \times P(R^d) \times P(R^d)) \times (D' \times P(R^d) \times P(R^d)) \to R$ be $\sigma((y_0, \bar{u}), (y_1, \overline{Mu})) = \log \frac{\|M(\theta_1)u(\theta_1)\|/\|u(\theta_1)\|}{\|M(\theta_0)u(\theta_0)\|/\|u(\theta_0)\|}$. For $\pi(\theta_0), \pi(\theta_1) \in P(R^d)$, denote $\sigma(W_0, W_0) = \log \frac{\|\mathbb{T}_0(\theta_1)\pi(\theta_1)\|/\|\pi(\theta_1)\|}{\|\mathbb{T}_0(\theta_0)\pi(\theta_0)\|/\|\pi(\theta_0)\|}$. Then

$$
\begin{aligned}
\mathbb{S}_n &= \log LR_n \\
&= \log \frac{\|M_n(\theta_1) \cdots M_1(\theta_1) M_0(\theta_1) \pi(\theta_1)\|}{\|M_n(\theta_0) \cdots M_1(\theta_0) M_0(\theta_0) \pi(\theta_0)\|} \\
&= \log \frac{\|\mathbb{T}_n(\theta_1)\pi(\theta_1)\|/\|\mathbb{T}_{n-1}(\theta_1)\pi(\theta_1)\|}{\|\mathbb{T}_n(\theta_0)\pi(\theta_0)\|/\|\mathbb{T}_{n-1}(\theta_0)\pi(\theta_0)\|} + \cdots \\
&\quad + \log \frac{\|\mathbb{T}_1(\theta_1)\pi(\theta_1)\|/\|\mathbb{T}_0(\theta_1)\pi(\theta_1)\|}{\|\mathbb{T}_1(\theta_0)\pi(\theta_0)\|/\|\mathbb{T}_0(\theta_0)\pi(\theta_0)\|} \\
&\quad + \log \frac{\|\mathbb{T}_0(\theta_1)\pi(\theta_1)\|/\|\pi(\theta_1)\|}{\|\mathbb{T}_0(\theta_0)\pi(\theta_0)\|/\|\pi(\theta_0)\|} \\
&= \sigma(W_{n-1}^\theta, W_n^\theta) + \cdots + \sigma(W_0^\theta, W_1^\theta) + \sigma(W_0^\theta, W_0^\theta)
\end{aligned}
$$
(2.15)

is an additive functional of the Markov chain $\{W_n^\theta, n \geq 0\}$.

**3. A nonlinear Markov renewal theory.** Note that $\{W_n^\theta, n \geq 0\}$ defined in (2.13) is a Markov chain on a general state space $D' \times P(R^d) \times P(R^d)$. In this section, abuse the notation a little bit and let $\{X_n, n \geq 0\}$ be a Markov chain on a general state space $\mathcal{X}$ with $\sigma$-algebra $\mathcal{A}$, which is irreducible with respect to a maximal irreducibility measure on $(\mathcal{X}, \mathcal{A})$ and is aperiodic. Let $S_n = \sum_{k=1}^n \xi_k$ be the additive component, taking values on the real line $R$, such that $\{(X_n, S_n), n \geq 0\}$ is a Markov chain on $\mathcal{X} \times R$ with transition

probability

$$P\{(X_{n+1}, S_{n+1}) \in A \times (B+s) | (X_n, S_n) = (x, s)\}$$
$$(3.1) \qquad = P\{(X_1, S_1) \in A \times B | (X_0, S_0) = (x, 0)\} = P(x, A \times B),$$

for all $x \in \mathcal{X}$, $A \in \mathcal{A}$ and $B \in \mathcal{B}(R)$ (:= Borel $\sigma$-algebra on $R$). The chain $\{(X_n, S_n), n \geq 0\}$ is called a *Markov random walk*. In this section, let $P_\nu$ ($E_\nu$) denote the probability (expectation) under the initial distribution on $X_0$ being $\nu$. If $\nu$ is degenerate at $x$, we shall simply write $P_x$ ($E_x$) instead of $P_\nu$ ($E_\nu$). We assume throughout this section that there exists a stationary probability distribution $\pi$, $\pi(A) = \int P(x, A) \, d\pi(x)$ for all $A \in \mathcal{A}$ and $E_\pi \xi_1 > 0$.

Let $\{Z_n = S_n + \eta_n, n \geq 0\}$ be a perturbed Markov random walk in the following sense: $S_n$ is a Markov random walk, $\eta_n$ is $\mathcal{F}_n$-measurable, where $\mathcal{F}_n$ is the $\sigma$-algebra generated by $\{(X_k, S_k), 0 \leq k \leq n\}$, and $\eta_n$ is *slowly changing*, that is, $\max_{1 \leq t \leq n} |\eta_t|/n \to 0$ in probability. Let $\{A = A(t; \lambda), \lambda \in \Lambda\}$ be a family of boundary functions for some index set $\Lambda$. Define

$$(3.2) \quad T = T_\lambda = \inf\{n \geq 1 : Z_n > A(n; \lambda)\}, \qquad \inf \varnothing = \infty \text{ for each } \lambda \in \Lambda.$$

It is easy to see that for all $\lambda > 0$, $T_\lambda < \infty$ with probability 1. This section concerns the approximations of the distribution of the overshoot and expected stopping time $E_\nu T$ as the boundary tends to infinity.

In the case of independent and identically distributed (i.i.d.) random variables $\xi_n$ with common positive mean, nonlinear renewal theory concerning boundary crossing times and its applications has been studied by Lai and Siegmund (1977, 1979), Woodroofe (1976, 1977) and Zhang (1988), among others. A good summary for this topic can be found in Woodroofe (1982) and Siegmund (1985) and references therein. For a perturbed Markov random walk with $E_\pi \xi_1 > 0$, Melfi (1992) generalized Lai and Siegmund's (1977) results to study the limiting distribution of the overshoot crossing a constant boundary. A multidimensional nonlinear first passage probability for perturbed Markov random walks can be found in Fuh and Lai (2001).

A Markov chain $\{X_n, n \geq 0\}$ on a state space $\mathcal{X}$ is called $V$-uniformly ergodic if there exists a measurable function $V : \mathcal{X} \to [1, \infty)$, with $\int V(x) \, d\pi(x) < \infty$, such that, for any Borel measurable function $h$ on $\mathcal{X}$ satisfying $\|h\|_V := \sup_x |h(x)|/V(x) < \infty$, we have

$$\lim_{n \to \infty} \sup_{x \in \mathcal{X}} \left\{ \frac{|E(h(X_n)|X_0 = x) - \int h(x) \, d\pi(x)|}{V(x)} : x \in \mathcal{X}, |h| \leq V \right\} = 0.$$

In this section we shall assume that $\{X_n, n \geq 0\}$ is $V$-uniformly ergodic. Under the irreducibility and aperiodicity assumption, $V$-uniform ergodicity implies that there exist $r > 0$ and $0 < \rho < 1$ such that for all $h$ and $n \geq 1$,

$$(3.3) \qquad \sup_{x \in \mathcal{X}} \frac{|E(h(X_n)|X_0 = x) - \int h(y) \, d\pi(y)|}{V(x)} \leq r \rho^n \|h\|_V;$$

see pages 382 and 383 of Meyn and Tweedie (1993). When $V \equiv 1$, this reduces to the classical uniform ergodicity condition.

The following assumptions for Markov chains will be used in this section:

A1. $\sup_x \{\frac{E(V(X_1))}{V(x)}\} < \infty$.
A2. $\sup_x E_x|\xi_1|^2 < \infty$ and $\sup_x \{\frac{E(|\xi_1|^r V(X_1))}{V(x)}\} < \infty$ for some $r \geq 1$.
A3. Let $\nu$ be an initial distribution of the Markov chain $\{X_n, n \geq 0\}$. Assume that for some $r \geq 1$,

$$(3.4) \qquad \sup_{\|h\|_V \leq 1} \left| \int_{x \in \mathcal{X}} h(x) E_x |\xi_1|^r \, d\nu(x) \right| < \infty.$$

A Markov random walk is called *lattice* with span $d > 0$ if $d$ is the maximal number for which there exists a measurable function $\gamma : \mathcal{X} \to [0, \infty)$ called the shift function, such that $P\{\xi_1 - \gamma(x) + \gamma(y) \in \{\ldots, -2d, -d, 0, d, 2d, \ldots\} | X_0 = x, X_1 = y\} = 1$ for almost all $x, y \in \mathcal{X}$. If no such $d$ exists, the Markov random walk is called *nonlattice*. A lattice random walk whose shift function $\gamma$ is identically 0 is called *arithmetic*.

To establish the nonlinear Markov renewal theorem, we shall make use of (3.1) in conjunction with the following extension of Cramér's (strongly nonlattice) condition [Götze and Hipp (1983), (2.5) on page 216]: There exists $\delta > 0$ such that for all $m, n = 1, 2, \ldots, \delta^{-1} < m < n$, and all $\theta \in R$ with $|\theta| \geq \delta$,

$$E_\pi | E\{\exp(i\theta(\xi_{n-m} + \cdots + \xi_{n+m})) | X_{n-m}, \ldots,$$
$$X_{n-1}, X_{n+1}, \ldots, X_{n+m}, X_{n+m+1}\} | \leq e^{-\delta}.$$

By using Markov renewal theory [Kesten (1974), Alsmeyer (1994), Fuh and Lai (2001) and Fuh (2004)] and Wald's equations for Markov random walks [Fuh and Lai (1998) and Fuh and Zhang (2000)], our approach is based on the investigation of the difference between $T_\lambda$ and a stopping time crossing linear boundaries with varying drift. That is, we first define

$$(3.5) \quad \tau := \tau(c, u) = \inf\{n \geq 1 : S_n - un > c\}, \qquad c \geq 0, \quad u \leq E_\pi \xi_1,$$

and establish the uniform integrability of $|T_\lambda - \tau(c_\lambda, d_\lambda)|^p$ for $p \geq 1$, for suitable $c_\lambda$ and $d_\lambda$. Then we derive nonlinear Markov renewal theory directly from parallel results in the linear case via the uniform integrabilities and the weak convergence of the overshoot.

Let $P_+^u(x, B \times R) = P_x\{X_{\tau(0,u)} \in B\}$ for $u \leq E_\pi \xi_1$, and denote the transition probability associated with the Markov random walk generated by the ascending ladder variable $S_{\tau(0,u)}$. Under the $V$-uniform ergodicity condition and $E_\pi \xi_1 > 0$, a similar argument as on pages 655–656 of Fuh and Lai (2001) yields that the transition probability $P_+^u(x, \cdot \times R)$ has an invariant measure $\pi_+^u$. Let $E_+^u$ denote expectation under $X_0$ having the initial distribution $\pi_+^u$. When $u = E_\pi \xi_1$, we denote $P_+^{E_\pi \xi_1}$ as $P_+$, and $\tau_+ = \tau(0, E_\pi \xi_1)$. Define

$$b = b_\lambda = \sup\{t \geq 1 : A(t, \lambda) \geq t E_\pi \xi_1\}, \qquad \sup \varnothing = 1, \tag{3.6}$$

$$d = d_\lambda = \left(\frac{\partial A}{\partial t}\right)(b_\lambda; \lambda), \tag{3.7}$$

$$\bar{d} = \sup\left\{\left(\frac{\partial A}{\partial t}\right)(t; \lambda); \ t \geq b_\lambda, \ \lambda \in \Lambda\right\}, \tag{3.8}$$

$$R = R_\lambda = Z_T - A(T; \lambda), \tag{3.9}$$

$$R(c, u) = S_{\tau(c,u)} - u\tau(c, u) - c, \qquad u \leq E_\pi \xi_1, \ c \geq 0, \tag{3.10}$$

$$r(u) = E_+^u R^2(0, u) / 2 E_+^u R(0, u), \qquad u \leq E_\pi \xi_1, \tag{3.11}$$

$$G(r, u) = \int_r^\infty P_+^u \{R(0, u) > s\} \, ds / E_+^u R(0, u), \qquad u \leq E_\pi \xi_1, \ r \geq 0. \tag{3.12}$$

We shall assume that $A(t; \lambda)$ is twice differentiable in $t$ and $b_\lambda$ is finite so that $d$ and $\bar{d}$ are well defined. The next theorem is a Blackwell-type nonlinear Markov renewal theorem. In the case of i.i.d. random variables, such a result has been developed by Lai and Siegmund (1977). Melfi (1992) has extended their result to the Markov case under a different ergodicity assumption as in this paper. Here, we consider a nonlinear boundary, extending Zhang's (1988) result under the $V$-uniform ergodicity assumption. Since for $1/2 < \alpha \leq 1$, $b^{-\alpha}(T - b) = o_{P_\nu}(1)$ implies (3.13) with $\gamma(b)/b^\alpha \to 0$, Theorem 1 implies Theorem 3 of Melfi (1992).

THEOREM 1. *Assume* A1 *holds, and* A2 *and* A3 *hold with* $r = 1$. *Let $\nu$ be an initial distribution on $X_0$. Suppose there exist functions $\rho(\delta) > 0$, $\sqrt{b} \leq \gamma(b) \leq b$, $\gamma(b)/b \to 0$ as $b \to \infty$, and a constant $d^* < E_\pi \xi_1 \in (0, \infty)$ such that*

$$(T_\lambda - b_\lambda)/\gamma(b_\lambda) = O_{P_\nu}(1) \qquad \text{as } b_\lambda \to \infty, \tag{3.13}$$

$$\lim_{n \to \infty} P_\nu \left\{\max_{1 \leq j \leq \rho(\delta)\gamma(n)} |\eta_{n+j} - \eta_n| \geq \delta\right\} = 0 \qquad \text{for any } \delta > 0, \tag{3.14}$$

$$\sup\left\{\left|\gamma^2(b)\left(\frac{\partial^2 A}{\partial t^2}\right)(t; \lambda)\right| : |t - b| \leq K\gamma(b), \ \lambda \in \Lambda\right\} < \infty \tag{3.15}$$

*for all $K > 0$*

*and*

$$\lim_{b_\lambda \to \infty} d_\lambda = d^*. \tag{3.16}$$

If $\xi_1 - d^*$ does not have an arithmetic distribution under $P_\nu$, then for any $r \geq 0$,

$$
\begin{aligned}
&P_\nu\{X_T \in B, R_\lambda > r\} \\
(3.17) \quad &= \frac{1}{E_+^{d^*} R(0, d^*)} \int_{x \in B} d\pi_+^{d^*}(x) \int_r^\infty P_+^{d^*}\{R(0, d^*) > s\} \, ds + o(1)
\end{aligned}
$$

as $b_\lambda \to \infty$.

In particular, $P_\nu\{R_\lambda > r\} = G(r, d^*) + o(1)$, as $b_\lambda \to \infty$ for any $r \geq 0$. If, in addition, $(T - b_\lambda)/\gamma(b_\lambda)$ converges in distribution to a random variable $W$ as $b_\lambda \to \infty$, then

$$(3.18) \quad \lim_{b_\lambda \to \infty} P_\nu\{R_\lambda > r, \ T_\lambda \geq b_\lambda + t\gamma(b_\lambda)\} = G(r, d^*) P_+^{d^*}\{W \geq t\},$$

for every real number $t$ with $P_+^{d^*}\{W = t\} = 0$.

The proof of Theorem 1 is given in Section 6.

To study uniform integrabilities of the powers of the differences for linear and nonlinear stopping times, we shall first give the regularity conditions on $\eta = \{\eta_n, n \geq 1\}$. The process $\eta$ is said to be *regular* with $p \geq 0$ and $1/2 < \alpha \leq 1$ if there exist a random variable $L$, a function $f(\cdot)$ and a sequence of random variables $U_n$, $n \geq 1$, such that

$$(3.19) \quad \eta_n = f(n) + U_n \quad \text{for } n \geq L \quad \text{and} \quad \sup_{x \in \mathcal{X}} E_x L^p < \infty,$$

$$(3.20) \quad \max_{1 \leq j \leq \sqrt{n}} |f(n+j) - f(n)| \leq K, \quad K < \infty,$$

$$(3.21) \quad \left\{ \max_{1 \leq j \leq n^\alpha} |U_{n+j}|^p, \ n \geq 1 \right\} \text{ is uniformly integrable,}$$

$$(3.22) \quad n^p \sup_{x \in \mathcal{X}} P_x \left\{ \max_{0 \leq j \leq n} U_{n+j} \geq \theta n^\alpha \right\} \to 0 \quad \text{as } n \to \infty, \text{ for all } \theta > 0,$$

and for some $w > 0$, $w < E_\pi \xi_1 - \bar{d}$ if $\alpha = 1$,

$$(3.23) \quad \sum_{n=1}^\infty n^{p-1} \sup_{x \in \mathcal{X}} P_x\{-U_n \geq wn^\alpha\} < \infty.$$

We shall set $f(n)$ to be the median of $\eta_n$ when $\eta$ is not regular and extend $f$ to a function on $[1, \infty)$ by linear interpolation. Therefore, we can define $\tau = \tau_\lambda = \tau(c_\lambda, d_\lambda)$ and $c_\lambda = b_\lambda(E_\pi \xi_1 - d_\lambda) - f(b_\lambda)$.

THEOREM 2. *Assume* A1 *holds, and* A2 *and* A3 *hold with* $r = p'(p + 1)/\alpha$ *for some* $p \geq 1$, $p' > 1$ *and* $1/2 < \alpha \leq 1$. *Suppose* $\eta$ *is regular with* $p \geq 1$, $1/2 < \alpha \leq 1$, and that there exist constants $\delta$ and $\mu^*$ with $0 < \delta < 1$ and $0 < \mu^* < E_\pi \xi_1$ such that

$$(3.24) \qquad b^p \sup_x P_x\{T \leq \delta b\} \to 0 \qquad \text{as } b \to \infty,$$

and

$$(3.25) \qquad \left(\frac{\partial A}{\partial t}\right)(t;\lambda) \leq \mu^*, \qquad t \geq \delta b, \ \lambda \in \Lambda.$$

(i) If $\sup_{x \in \mathcal{X}} E_x\{|\xi_1|^{2pp'}\} < \infty$ for some $p' > 1$ and for any $K > 0$,

$$(3.26) \quad \sup\left\{\left|b_\lambda\left(\frac{\partial^2 A}{\partial t^2}\right)(t;\lambda)\right| : b_\lambda - Kb_\lambda^\alpha \leq t \leq b_\lambda + Kb_\lambda^\alpha, \ \lambda \in \Lambda\right\} < \infty,$$

then

$$(3.27) \qquad \{|T_\lambda - \tau_\lambda|^p; \ \lambda \in \Lambda\} \text{ is uniformly integrable under } P_\nu.$$

(ii) If $\partial^2 A/\partial t^2 = 0$, then (3.27) still holds without the condition $\sup_x E_x\{|\xi_1|^{2pp'}\} < \infty$.

The proof of Theorem 2 is given in Section 6.

We need the following notation and definitions before stating Theorem 3. For a given Markov random walk $\{(X_n, S_n), n \geq 0\}$, let $\nu$ be an initial distribution of $X_0$ and define $\nu^*(B) = \sum_{n=0}^\infty P_\nu(X_n \in B)$ on $\mathcal{A}$. Let $g = E(\xi_1|X_0, X_1)$ and $E_\pi|g| < \infty$. Define operators $\mathbf{P}$ and $\mathbf{P}_\pi$ by $(\mathbf{P}g)(x) = E_x g(x, X_1, \xi_1)$ and $\mathbf{P}_\pi g = E_\pi g(X_0, X_1, \xi_1)$, respectively, and set $\bar{g} = \mathbf{P}g$. We shall consider solutions $\Delta(x) = \Delta(x;g)$ of the Poisson equation

$$(3.28) \qquad (I - \mathbf{P})\Delta = (I - \mathbf{P}_\pi)\bar{g}, \qquad \nu^*\text{-a.s.}, \ \mathbf{P}_\pi \Delta = 0,$$

where $I$ is the identity operator. Under conditions A1–A4, it is known [Theorem 17.4.2 of Meyn and Tweedie (1993)] that the solution $\Delta$ of (3.28) exists and is bounded.

THEOREM 3. *Assume A1 holds, and A2 and A3 hold with $r = 2 + p$ for some $p > 1$. Let $\nu$ be an initial distribution such that $E_\nu V(X_0) < \infty$. Suppose that*

$$(3.29) \quad \lim_{n \to \infty} \sup_{x \in \mathcal{X}} P_x\left\{\max_{1 \leq j \leq \sqrt{n}} |\eta_{n+j} - \eta_j| \geq \delta\right\} = 0 \qquad \text{for any } \delta > 0,$$

$$(3.30) \qquad \eta_n = f(n) + U_n \qquad \text{for any } n \geq L,$$

*and that there exist constants $d_1^* < E_\pi \xi_1$ and $d_2^*$ such that*

$$(3.31) \qquad \lim_{n \to \infty} \max_{0 \leq j \leq \sqrt{n}} |f(n+j) - f(n)| = 0,$$

$(3.32)$ $U_n$ *converges in distribution to an integrable random variable $U$,*

$$(3.33) \qquad \lim_{b_\lambda \to \infty} d_\lambda = d_1^* \quad \text{and} \quad \xi_1 - d_1^* \text{ is nonarithmetic under } P_\nu,$$

and for any constant $K > 0$,

$$(3.34) \quad \lim_{b_\lambda \to \infty} \sup\left\{\left|b_\lambda\left(\frac{\partial^2 A}{\partial t^2}\right)(t;\lambda) - d_2^*\right| : (t - b_\lambda)^2 \leq Kb_\lambda\right\} = 0.$$

If $\{|T_\lambda - \tau_\lambda|; \lambda \in \Lambda\}$ is uniformly integrable, then

$$(3.35) \quad E_\nu T_\lambda = b_\lambda - (E_\pi \xi_1 - d_\lambda)^{-1} f(b_\lambda) + C_0 + o(1) \quad \text{as } b_\lambda \to \infty,$$

where

$$(3.36) \quad C_0 = (E_\pi \xi_1 - d_1^*)^{-1}\bigg(r(d_1^*) + (E_\pi \xi_1 - d_1^*)^{-2} d_2^* \sigma^2/2 - E_\pi U$$
$$- \int \Delta(x)\, d(\pi_+^{d^*}(x) - \nu(x))\bigg).$$

The proof of Theorem 3 is given in Section 6.

When $A(t, \lambda) = \lambda$, we have the following:

COROLLARY 1. *Under the assumptions of Theorem 3, as $\lambda \to \infty$,*

$$(3.37) \quad E_\nu T_\lambda = (E_\pi \xi_1)^{-1}\bigg(\lambda + E_{\pi_+} S_{\tau_+}^2/2S_{\tau_+} - f(\lambda/E_\pi \xi_1) - E_\pi U$$
$$- \int \Delta(x)\, d(\pi_+(x) - \nu(x))\bigg) + o(1).$$

**4. Asymptotic optimality of the SRP detection procedure.** For ease of notation, let $\mathcal{X} := D' \times P(R^d) \times P(R^d)$ be the state space of the Markov chain $\{W_n^\theta, n \geq 0\}$ defined in (2.13). Denote $w := (y, \bar{u}, \bar{u})$ and $\tilde{w} := (y_0, \pi, \pi)$, where $y_0 = (x_0, \pi) \in D'$ and $x_0$ is the initial state of $X_0$ taken from $\pi$. To prove the asymptotic optimality of the SRP rule in hidden Markov models, the following condition C will be assumed throughout this paper.

C1. For each $\theta \in \Theta$, the Markov chain $\mathbf{X} = \{X_n, n \geq 0\}$ is ergodic (positive recurrent, irreducible and aperiodic) on a finite state space $D = \{1,\ldots,d\}$. Moreover, the Markov chain $\{(X_n, \xi_n), n \geq 0\}$ is irreducible, aperiodic and $V$-uniformly ergodic for some $V$ on $D'$ with A1 and A2 holding. We also assume the Markov chain $\{X_n, n \geq 0\}$ has stationary probability $\Gamma$ with probability density $\pi_x(\theta) f(\cdot; \varphi_x(\theta))$ with respect to $\mu$.

C2. For each $\theta \in \Theta$, the random matrices $M_0(\theta)$ and $M_1(\theta)$ defined in (2.9) and (2.10) are invertible $\mathbb{P}^\theta$ almost surely and

$$\sup_{(x,\xi_0) \in D \times R} E_x^\theta \left|\sum_{x,y=1}^d \pi_x(\theta) f(\xi_0; \varphi_x(\theta)) p_{xy}(\theta) \xi_1 f(\xi_1; \varphi_y(\theta)|\xi_0)\right| < \infty.$$

The construction of the SRP rule and the proof of its asymptotic optimality can be split into two steps. We first prove that it is a limit of Bayes rules, and then we prove the asymptotic optimality. To this end, let us consider the Bayesian formulation of change point detection in a hidden Markov model and denote it by $B(\beta, p, c, \tilde{w})$. That is, we assume the initial state of $W_0$ is $\tilde{w}$ and suppose $\omega$ has a prior distribution

$$\mathbb{P}_{\tilde{w}}(\omega = 0) = \beta \quad \text{and} \quad \mathbb{P}_{\tilde{w}}(\omega = n) = (1-\beta)p(1-p)^{n-1} \quad \text{for } n \geq 1,$$

where $p$ and $\beta$ are known constants with $0 < p \leq 1$, $0 \leq \beta \leq 1$. The parameter $\omega$ is the (unknown) point of change of the process from a hidden Markov model.

Let $N$ be a stopping time adapted to the system of $\sigma$-algebras $\{\mathcal{F}_n\}_{n=0}^\infty$, where $\mathcal{F}_0$ is the natural $\sigma$-algebra $\{\varnothing, \mathcal{X}\}$ and $\mathcal{F}_n = \sigma(\mathcal{F}_0, W_0, W_1, \ldots, W_n)$. Following the formulation of Shiryayev (1963, 1978) and its modification given by Yakir (1994) for a finite state Markov chain, the risk associated with the detection policy $N$ is

$$(4.1) \qquad \rho(N, \omega) = \mathbb{P}_{\tilde{w}}(N < \omega) + c\,\mathbb{E}_{\tilde{w}}(N - \omega)^+,$$

where $a^+$ denotes $\max\{a, 0\}$, and $c > 0$ is a fixed constant.

DEFINITION 1. For a given pair $(p, \tilde{w}) \in (0, 1] \times \mathcal{X}$, we call a stopping time $N^*$ a $B(\beta, p, c, \tilde{w})$-Bayes time if

$$\rho(N^*, \omega) = \inf \rho(N, \omega),$$

where inf is taken over the class of all proper stopping times.

The following proposition characterizes the structure of the Bayes rule in hidden Markov models. Since the proof of the proposition is similar to Shiryayev's proof, it is omitted.

PROPOSITION 1. *Let* $0 < p \leq 1$, $c > 0$ *and let*

$$\delta_n = \delta_n(p, \tilde{w}) = \mathbb{P}_{\tilde{w}}(\omega \leq n | \mathcal{F}_n)$$

*be the posterior probability that the next observation is governed by* $P^{\theta_1}$. *There exists a function* $A_p(\cdot)$, *defined on* $\mathcal{X}$, *such that the stopping time*

$$(4.2) \qquad N_{A,p} = \inf\{n \geq 0 : \delta_n(p, \tilde{w}) \geq A_p(W_n)\}$$

*is the* $B(\beta, p, c, \tilde{w})$-*Bayes rule. Moreover,* $A_p(\cdot)$ *does not depend on* $\beta$ *or on* $\tilde{w}$.

REMARK 1. Proposition 1 remains correct when the initial pair $(p, \tilde{w})$ is random (according to a measure $\phi$). Again, the threshold function does not depend on the initial state. (Notice that the stopping time does depend on the distribution of the initial state through the dependence on the initial state of the probability of a change.) The structure of the Bayes rule plays a crucial role in the development of the optimal detection time in the non-Bayesian setting.

Denote

$$r(x) = \frac{x}{(1-x)p}, \qquad q = 1 - p, \tag{4.3}$$

and let

$$LR_{n,p} = r(\beta)\frac{LR_n}{q} + \sum_{k=0}^{n} \frac{LR_n^k}{q}. \tag{4.4}$$

It is convenient to reformulate the stopping time $N_{A,p}$ in terms of a different sequence of statistics. By using the same idea as Lemma 2 of Pollak (1985), it follows that

$$\delta_n(p, \tilde{w}) = \frac{LR_{n,p}}{LR_{n,p} + 1/p}. \tag{4.5}$$

Since the function $y/(y + 1/p)$ is a monotone function in $y$, the Bayesian stopping time can be rewritten in terms of $LR_{n,p}$,

$$\begin{aligned}
N_{A,p} &= \inf\{n \geq 0 : \delta_n(p, \tilde{w}) \geq A_p(W_n)\} \\
&= \inf\{n \geq 0 : LR_{n,p} \geq B_p(W_n)\} = N_{B,p},
\end{aligned} \tag{4.6}$$

where $B_p(\cdot) = r(A_p(\cdot))$. For consistency of notation, we will use $N_{B,p}$ instead of $N_{A,p}$ in the sequel.

THEOREM 4. *Assume* C1 *and* C2 *hold. Suppose that the $P_\infty$-distribution of $LR_1$ is nonarithmetic.*

(i) *There exists a $\delta > 0$ such that for any $\delta < b = \log B < \infty$, there exist a constant $0 < c^* < \infty$ and a sequence $\{p_i, c_i\}_{i=1}^\infty$ with $p_i \to 0$, $c_i \to c^*$ as $i \to \infty$ such that the stopping time $N_b$ defined in (2.7) is a limit as $i \to \infty$ of Bayes rules for $B$ $(\beta = 0, p = p_i, c = c_i, \tilde{w})$.*

(ii) *For any set of Bayes problems $B(\beta, p, c, \tilde{w})$ with $\beta = 0$, $p \to 0$, $c \to c^*$,*

$$\limsup_{p \to 0, \ c \to c^*} \frac{1 - E\rho(N_{B,p}, \omega)}{1 - E\rho(N_b, \omega)} = 1, \tag{4.7}$$

*where the expectation is taken in the Bayes problems $B(0, p, c, \tilde{w})$.*

(iii) *For any $1 \leq \gamma < \infty$, there exists a unique $1 < b = \log B < \infty$ such that $\gamma = E_\infty N_b$.*

The proof of Theorem 4 is given in Section 7.

After understanding the structure of the Bayes rules for detecting a change in hidden Markov models and the characteristics of the limits of such rules, we can turn our attention to the problem of detecting a change in a non-Bayesian setting. To study randomization of the initial for the SRP change point detection rule, we need the following notation first.

For $0 \leq k \leq n$, let

$$(4.8) \qquad R_{n,p} := \sum_{k=0}^{n} \frac{1}{q} \frac{p_n(\xi_k, \xi_{k+1}, \ldots, \xi_n; \theta_1)}{p_n(\xi_k, \xi_{k+1}, \ldots, \xi_n; \theta_0)}.$$

Note that $R_{n,p} = LR_{n,p}$ when $\beta = 0$. By using the same notation as that in Section 2, for $y_0, y_1 \in D'$, $\bar{u} = \overline{u(\theta)} = (\overline{u(\theta_0)}, \overline{u(\theta_1)}) \in P(R^d) \times P(R^d)$ and $M = M(y_0, y_1) = M(\theta) = (M(\theta_0), M(\theta_1)) \in Gl(d, R) \times Gl(d, R)$, let $\beta \colon (D' \times P(R^d) \times P(R^d)) \times (D' \times P(R^d) \times P(R^d)) \to R$ be $\beta((y_0, \bar{u}), (y_1, \overline{Mu})) = \frac{\|M(\theta_1)u(\theta_1)\|/\|u(\theta_1)\|}{\|M(\theta_0)u(\theta_0)\|/\|u(\theta_0)\|}$. For $\pi(\theta_0), \pi(\theta_1) \in P(R^d)$, denote $\beta(W_0, W_0) = \frac{\|\mathbb{T}_0(\theta_1)\pi(\theta_1)\|/\|\pi(\theta_1)\|}{\|\mathbb{T}_0(\theta_0)\pi(\theta_0)\|/\|\pi(\theta_0)\|}$. Then

$$\begin{aligned}
&\frac{p_n(\xi_0, \xi_1, \ldots, \xi_n; \theta_1)}{p_n(\xi_0, \xi_1, \ldots, \xi_n; \theta_0)} \\
&= \frac{\|M_n(\theta_1) \cdots M_1(\theta_1) M_0(\theta_1) \pi(\theta_1)\|}{\|M_n(\theta_0) \cdots M_1(\theta_0) M_0(\theta_0) \pi(\theta_0)\|} \\
(4.9) \quad &= \frac{\|\mathbb{T}_n(\theta_1)\pi(\theta_1)\|/\|\mathbb{T}_{n-1}(\theta_1)\pi(\theta_1)\|}{\|\mathbb{T}_n(\theta_0)\pi(\theta_0)\|/\|\mathbb{T}_{n-1}(\theta_0)\pi(\theta_0)\|} \cdots \\
&\quad \times \frac{\|\mathbb{T}_1(\theta_1)\pi(\theta_1)\|/\|\mathbb{T}_0(\theta_1)\pi(\theta_1)\|}{\|\mathbb{T}_1(\theta_0)\pi(\theta_0)\|/\|\mathbb{T}_0(\theta_0)\pi(\theta_0)\|} \cdot \frac{\|\mathbb{T}_0(\theta_1)\pi(\theta_1)\|/\|\pi(\theta_1)\|}{\|\mathbb{T}_0(\theta_0)\pi(\theta_0)\|/\|\pi(\theta_0)\|} \\
&= \beta(W_{n-1}^\theta, W_n^\theta) \cdots \beta(W_0^\theta, W_1^\theta) \cdot \beta(W_0^\theta, W_0^\theta)
\end{aligned}$$

is a product of the functional for the Markov chain $\{W_n^\theta, n \geq 0\}$. Therefore, (4.8) can be rewritten as

$$(4.10) \quad R_{n,p} = \sum_{k=0}^{n} \frac{1}{q} \beta(W_{n-1}^\theta, W_n^\theta) \cdots \beta(W_{k-1}^\theta, W_k^\theta) \qquad \text{where } W_{-1}^\theta = W_0^\theta.$$

Define

$$R_{n+1,p} = \beta(W_{n-1}^\theta, W_n^\theta) \frac{1}{q}(1 + R_{n,p}), \qquad N_{q,b} := \inf\{n : R_{n,p} \geq B\},$$

$$F_n(s, w) = \mathbb{P}_\infty(R_{n+1,p} \leq s \mid N_{q,b} > n, W_n = w),$$

$$\rho(t, s, w, w') = \mathbb{P}_\infty(R_{n+1,p} \leq s, W_{n+1} \in dw' \mid R_{n,p} = t, N_{q,b} > n+1, W_n = w),$$

$$\zeta(t, w, w') = \mathbb{P}_\infty(N_{q,b} > n+1, W_{n+1} \in dw' \mid R_{n,p} = t, N_{q,b} > n, W_n = w).$$

For a given set of nonnegative boundary points $B = \{B(w) : w \in \mathcal{X}\}$ (infinity is not excluded), consider the set $S_B = \{(r, w) : w \in \mathcal{X}, 0 < r < B(w)\}$. Let $\mathcal{F}_B$ be the set of distribution functions with support in $S_B$. Let $T_B$ be the transformation on $\mathcal{F}_B$ defined by

$$(4.11) \quad T_B F(r, w) = \frac{1}{Q(F)} \int_{w' \in \mathcal{X}} \int_0^{B(w')} \rho(t, r, w, w') \zeta(t, w, w') \times \mathbb{P}(w, dw') \, dF(t, w'),$$

where

$$(4.12) \quad Q(F) = \int_{w, w' \in \mathcal{X}} \int_0^{B(w')} \zeta(t, w, w') \mathbb{P}(w, dw') \, dF(t, w').$$

The idea behind (4.11) and (4.12) comes from iterated random functions, which Pollak (1985) used to define a change point detection rule in the independent case. Here $F_n(s, w)$ is driven by the Markov chain $\{W_n^\theta, n \geq 0\}$, and, hence, in the domain of Markovian iterated random functions. Under some regularity conditions on the Markov chain $\{W_n^\theta, n \geq 0\}$, and the continuity property for the iterated random functions, we will show in Lemma 8 that for each $B$ there is an associated set of invariant measures $\Phi_B$, that is, $T_B \phi = \phi$ for all $\phi \in \Phi_B$. Let $p = 1 - q$ and define $\tilde{\phi}$ as

$$d\tilde{\phi}(s, w) = \frac{(1 + ps) \, d\phi(s, w)}{\int_{w \in \mathcal{X}} \int_0^{B(w')} (1 + pt) \, d\phi(t, w)}.$$

It is easy to see that if the distribution of $R_{0,p}$ is $\tilde{\phi}$, then the distribution of $R_{0,p}$ conditional on $\{\omega > 0\}$ is $\phi$. Note that $\phi$ depends on $p$. By using the same argument as that in Theorem 4, we can choose a subsequence $\{T_B^i, p_i, c_i, \phi_i\}$ such that as $i \to \infty$, $p_i \to 0$, $c_i \to c^*$ and $\phi_i$ converges in distribution to a limit $\psi$.

Given the value of the initial state $W_0 = \tilde{w}$, the initial $(R_0^*, \tilde{w})$ is simulated from the distribution $\psi$, conditioned on the event $\{W_0 = \tilde{w}\}$. Define recursively

$$(4.13) \quad R_{n+1}^* = \beta(W_{n-1}^\theta, W_n^\theta)(1 + R_n^*).$$

Denote $b = \log B$, and define the SRP rule

$$(4.14) \quad N_b^\psi := \inf\{n : R_n^* \geq B\} = \inf\{n : \log R_n^* \geq b\}.$$

Notice that each one of these detection policies is an "equalizer rule" in the sense that

$$(4.15) \quad \mathbb{E}_k(N_b^\psi - k + 1 | N_b^\psi \geq k - 1) = \mathbb{E}_1 N_b^\psi,$$

for all $k > 1$. The same is true for the case where $\psi$ has atoms on the boundary, since the randomization law is time independent.

Note that the threshold of the Bayes rule (4.6) depends on the current state of the Markov chain, while the threshold of the SRP rule (4.14) is a constant. We claim in Lemma 7 and Theorem 5 that the difference between these two rules is $o(1)$ as $\gamma \to \infty$, by which we prove the conjecture raised by Yakir (1994) for finite state Markov chains.

THEOREM 5. *Assume* C1 *and* C2 *hold. Suppose that the* $P_\infty$-*distribution of* $LR_1$ *is nonarithmetic. Then for any* $1 < \gamma < \infty$, *there exist a constant* $\delta < b = \log B < \infty$ *and a probability measure* $\psi$ *such that* $\gamma = E_\infty N_b^\psi$ *and such that if* $N$ *is any stopping time which satisfies* $E_\infty N \geq \gamma$, *then*

$$(4.16) \quad \sup_{1 \leq \omega < \infty} E_\omega(N - \omega | N \geq \omega) \geq \sup_{1 \leq \omega < \infty} E_\omega(N_b^\psi - \omega | N_b^\psi \geq \omega) + o(1),$$

*where* $o(1) \to 0$ *as* $\gamma \to \infty$, $E_\omega(N_b^\psi - \omega | N_b^\psi \geq \omega)$ *is a constant for* $1 < \omega < \infty$.

The proof of Theorem 5 is given in Section 7.

**5. Asymptotic approximations for the average run length.** Since $N_b^\psi$ is an equalizer rule in the sense of $\mathbb{E}_k(N_b^\psi - k + 1 | N_b^\psi \geq k - 1) = \mathbb{E}_1 N_b^\psi$ by (4.15), in this section we consider only the approximation of $\mathbb{E}_1 N_b^\psi$. For $\theta = \theta^0$ or $\theta^1$, let $\pi^\theta$ denote the stationary distribution of $\{X_n, n \geq 0\}$ under $P^\theta$. For given $\mathbb{P}^{\theta_0}$ and $\mathbb{P}^{\theta_1}$, define the Kullback–Leibler information numbers as (4.2) of Fuh (2003),

$$(5.1) \quad \begin{aligned} K(\mathbb{P}^{\theta_0}, \mathbb{P}^{\theta_1}) &= \mathbb{E}_{\mathbb{P}^{\theta_0}} \left( \log \frac{\|M_1(\theta_0) M_0(\theta_0) \pi^{\theta_0}\|}{\|M_1(\theta_1) M_0(\theta_1) \pi^{\theta_1}\|} \right), \\ K(\mathbb{P}^{\theta_1}, \mathbb{P}^{\theta_0}) &= \mathbb{E}_{\mathbb{P}^{\theta_1}} \left( \log \frac{\|M_1(\theta_1) M_0(\theta_1) \pi^{\theta_1}\|}{\|M_1(\theta_0) M_0(\theta_0) \pi^{\theta_0}\|} \right), \end{aligned}$$

where $\mathbb{P}^{\theta_0}$ ($\mathbb{P}^{\theta_1}$) denotes the probability of the Markov chain $\{W_n^{\theta_0}, n \geq 0\}$ ($\{W_n^{\theta_1}, n \geq 0\}$), and $\mathbb{E}_{\mathbb{P}^{\theta_0}}$ ($\mathbb{E}_{\mathbb{P}^{\theta_1}}$) refers to the expectation under $\mathbb{P}^{\theta_0}$ ($\mathbb{P}^{\theta_1}$).

In the rest of this section we will impose the following mild condition on the Kullback–Leibler information numbers:

$$(5.2) \quad 0 < K(\mathbb{P}^{\theta_0}, \mathbb{P}^{\theta_1}) < \infty \quad \text{and} \quad 0 < K(\mathbb{P}^{\theta_1}, \mathbb{P}^{\theta_0}) < \infty.$$

To derive a second-order approximation for the average run length of the SRP rule, we will apply relevant results from nonlinear Markov renewal theory developed in Section 3. For this purpose, we rewrite the stopping time $N_b := N_b^\psi$ (we delete $\psi$ in this section for simplicity) in the form of a Markov random walk crossing a constant threshold plus a nonlinear term that is slowly changing. Note that the stopping time $N_b$ can be written in the form

$$(5.3) \quad N_b = \inf\{n \geq 1 : \mathbb{S}_n + \eta_n \geq b\}, \qquad b = \log B,$$

where $\mathbb{S}_n$ is a Markov random walk defined in (2.15) with mean $\mathbb{E}_1 \mathbb{S}_1 = K(\mathbb{P}^{\theta_1}, \mathbb{P}^{\theta_0})$, and

$$\eta_n = \log\left\{1 + \sum_{k=1}^{n-1} e^{-\mathbb{S}_k}\right\}. \tag{5.4}$$

For $b > 0$, define

$$N_b^* = \inf\{n \geq 1 : \mathbb{S}_n \geq b\}, \tag{5.5}$$

and let $R_b = \mathbb{S}_{N_b^*} - b$ (on $\{N_b^* < \infty\}$) denote the overshoot of the statistic $\mathbb{S}_n$ crossing the threshold $b$ at time $n = N_b^*$. When $b = 0$, we denote $N_b^*$ in (5.5) as $N_+^*$. For given $\tilde{w} \in \mathcal{X}$, let

$$G(y) = \lim_{b \to \infty} \mathbb{P}_1\{R_b \leq y | W_0 = \tilde{w}\} \tag{5.6}$$

be the limiting distribution of the overshoot. It is known [cf. Theorem 1 of Fuh (2004)] that

$$\lim_{b \to \infty} \mathbb{E}_1(R_b | W_0 = \tilde{w}) = \int_0^\infty y \, dG(y) = \frac{\mathbb{E}_{m_+} S_{N_+^*}^2}{2\mathbb{E}_{m_+} S_{N_+^*}},$$

where $m_+$ is defined in the same way as $\pi_+$ defined in the paragraph before (3.6) in Section 3.

Note that by (5.3),

$$\mathbb{S}_{N_b} = b - \eta_{N_b} + \chi_b \qquad \text{on } \{N_b < \infty\},$$

where $\chi_b = \mathbb{S}_{N_b} + \eta_{N_b} - b$ is the overshoot of $\mathbb{S}_n + \eta_n$ crossing the boundary $b$ at time $N_b$. Taking the expectations on both sides, and applying Wald's identity for products of Markovian random matrices [Theorem 2 of Fuh (2003)], we obtain

$$\begin{aligned}K(\mathbb{P}^{\theta_1}, \mathbb{P}^{\theta_0})\mathbb{E}_1(N_b | W_0 = \tilde{w}) - \int_{\mathcal{X}} \Delta(w') \, dm_+(w') + \Delta(\tilde{w}) \\ = \mathbb{E}_1(\mathbb{S}_{N_b} | W_0 = \tilde{w}) = b - \mathbb{E}_1(\eta_{N_b} | W_0 = \tilde{w}) + \mathbb{E}_1(\chi_b | W_0 = \tilde{w}),\end{aligned} \tag{5.7}$$

where $\Delta : \mathcal{X} \to R$ solves the Poisson equation

$$\mathbb{E}_w \Delta(W_1) - \Delta(w) = \mathbb{E}_w \mathbb{S}_1 - \mathbb{E}_m \mathbb{S}_1 \tag{5.8}$$

for almost every $w \in \mathcal{X}$ with $\mathbb{E}_m \Delta(W_1) = 0$.

The crucial observations are that the sequence $\{\eta_n, n \geq 1\}$ is slowly changing, and that $\eta_n$ converges $\mathbb{P}_1$-a.s. as $n \to \infty$ to the random variable

$$\eta = \log\left\{1 + \sum_{k=1}^{\infty} e^{-\mathbb{S}_k}\right\} \tag{5.9}$$

with finite expectation $\mathbb{E}_{m_+}\eta$. Here the expectation $\mathbb{E}_{m_+}$ is taken under $\omega = 1$ and the initial distribution of $W_0$ is $m_+$; we omit 1 for simplicity. An important consequence of the slowly changing property is that, under mild conditions, the limiting distribution of the overshoot of a Markov random walk over a fixed threshold does not change by the addition of a slowly changing nonlinear term (see Theorem 1).

THEOREM 6. *Assume* C1 *and* C2 *hold. Let* $\xi_0, \xi_1, \ldots, \xi_n$ *be a sequence of random variables from a hidden Markov model* $\{\xi_n, n \geq 0\}$. *Assume that* $\mathbb{S}_1$ *is nonarithmetic with respect to* $\mathbb{P}_\infty$ *and* $\mathbb{P}_1$. *If* $0 < K(\mathbb{P}^{\theta_1}, \mathbb{P}^{\theta_0}) < \infty$, $0 < K(\mathbb{P}^{\theta_0}, \mathbb{P}^{\theta_1}) < \infty$, *and* $\mathbb{E}_1|\mathbb{S}_1|^2 < \infty$, *then for* $\tilde{w} \in \mathcal{X}$, *as* $b \to \infty$,

$$
\begin{aligned}
&\mathbb{E}_1(N_b | W_0 = \tilde{w}) \\
&\quad = \frac{1}{K(\mathbb{P}^{\theta_1}, \mathbb{P}^{\theta_0})} \\
&\qquad \times \left( b - \mathbb{E}_{m_+}\eta + \frac{\mathbb{E}_{m_+} S^2_{N^*_+}}{2\mathbb{E}_{m_+} S_{N^*_+}} - \int_{\mathcal{X}} \Delta(w)\, dm_+(w) + \Delta(\tilde{w}) \right) + o(1).
\end{aligned}
$$
(5.10)

The proof of Theorem 7 is given in Section 8.

REMARK 2. The constants $\mathbb{E}_{m_+} S^2_{N^*_+}/2\mathbb{E}_{m_+} S_{N^*_+}$ and $\mathbb{E}_{m_+}\eta$ are the subject of the nonlinear renewal theory. The constant $-\int_{\mathcal{X}} \Delta(w)\, dm_+(w) + \Delta(\tilde{w})$ is due to Markovian dependence via the Poisson equation (5.8). Obviously, this bound is asymptotically accurate when $K(\mathbb{P}^{\theta_1}, \mathbb{P}^{\theta_0}) \to 0$.

**6. Proofs of Theorems 1–3.** We will use the same notation as that in Section 3 unless specifically mentioned.

PROOF OF THEOREM 1. To prove Theorem 1, we can make use of Theorem 3.1 for the one-dimensional case in Fuh and Lai (2001) as in the case of i.i.d. $\xi_n$ [see Theorem 1 of Zhang (1988)]. The details are omitted. □

To prove Theorem 2, we need some lemmas first.

LEMMA 1. *Let* $\nu$ *be an initial distribution such that* $E_\nu V(X_0) < \infty$. *Let* $\tau(c, u)$ *be defined by* (3.5), *and let* $p$ *and* $\alpha$ *be two constants with* $p \geq 1$ *and* $1/2 < \alpha \leq 1$.

(i) *If* $E_\nu(|\xi_1|^{p'(p+1)/\alpha} V(X_0)) < \infty$ *for some* $p' > 1$, *then*

$$
\sum_{n=1}^{\infty} n^{p-1} P_\nu \left\{ \max_{j \leq n} |S_j - j| \geq \gamma n^\alpha \right\} < \infty \quad \text{for any } \gamma > 0.
$$
(6.1)

(ii) If $E_\nu(|\xi_1|^{2pp'}V(X_0)) < \infty$ for some $p' > 1$, then for any $K > 0$,

$$\{((\tau(c,u) - (1-u)^{-1}c)^2/c)^p;\ c \geq 1,\ K^{-1} \leq 1-u \leq K\}$$

is uniformly integrable under $P_\nu$.

PROOF. By Theorem 16.0.1 of Meyn and Tweedie (1993), the $V$-uniform ergodicity condition is equivalent to the fact that there exist an extended real-valued function $w: \mathcal{X} \to [1, \infty)$, a measurable set $C$ and constants $\gamma > 0$, $b < \infty$, such that

$$\int_{\mathcal{X}} w(y)P(x,dy) - w(x) \leq -\gamma w(x) + bI_C(x) \qquad \text{for } x \in \mathcal{X},$$

where $w$ is equivalent to $V$ in the sense that for some $c \geq 1$, $c^{-1}V \leq w \leq cV$. Denote $\Delta_1 = \Delta$ as defined in (3.28). Let $\tilde{g} = E(d_1^2|X_0, X_1)$ and $\Delta_2(x; g) = \Delta(x; \tilde{g})$. Since $\int V(x)\,d\pi(x) < \infty$, A2 implies that there exists $0 < c < \infty$ such that for all $x \in \mathcal{X}$, $E(|\xi_1|^2|X_0 = x) < cV(x)$. By Theorem 17.4.2 of Meyn and Tweedie (1993), the solution $\Delta_r$ satisfies $\Delta_r \leq R_r(V(x) + 1)$ for $r = 1, 2$. This implies that $\sup_i E_\nu|\Delta_r(X_i; g)| < R_r \sup_i E_\nu(V(X_i) + 1) < \infty$ for $r = 1, 2$. Therefore, the conditions of Theorem 2 of Fuh and Zhang (2000) hold, and, hence, the quick convergence (i) follows from Theorem 2 of Fuh and Zhang (2000).

(ii) The proof of (ii) can be derived from (i) easily. □

Following Lemmas 2 and 3 in Zhang (1988), we have the following:

LEMMA 2. *Suppose that $\eta$ is regular with $p \geq 1$ and $1/2 < \alpha \leq 1$ and that conditions (3.24) and (3.25) hold. If $E_\nu(|\xi_1|^{p'(p+1)/\alpha}V(X_0)) < \infty$ for some $p' > 1$, then*

$$\lim_{b \to \infty} b^p P_\nu\{T_\lambda \leq b - \gamma b^\alpha\} = 0 \qquad \text{for any } \gamma > 0.$$

LEMMA 3. *Suppose that $\eta$ is regular with $p \geq 1$ and $1/2 < \alpha \leq 1$ and that condition (3.25) holds. Denote $n^* = [b + Kb^\alpha]$. If $E_\nu(|\xi_1|^{p'(p+1)/\alpha}V(X_0)) < \infty$ for some $p' > 1$, then there exists a constant $K > 0$ such that*

$$\lim_{b \to \infty} \sum_{n=n^*}^{\infty} n^{p-1} P_\nu\{T_\lambda > n\} = 0.$$

LEMMA 4. *Under the conditions of Theorem 2(i),*

$$\{((T_\lambda - \tau)^+)^p; \lambda \in \Lambda\} \text{ is uniformly integrable under } P_\nu.$$

PROOF. For ease of notation, let $T = T_\lambda$, $n_1 = [b - \gamma b^\alpha]$, $n^* = [b + Kb^\alpha]$, $T' = \max(n_1, \min(T, n^*))$ and $\tau' = \max(n_1, \min(\tau, n^*))$. By Lemmas 1(i), 2 and 3,

$$\lim_{b\to\infty} E_\pi |T - T'|^p = \lim_{b\to\infty} E_\pi |\tau - \tau'|^p = 0. \tag{6.2}$$

Let $\gamma' = (1 - \mu^*)/5$, where $\mu^*$ is defined in Theorem 2. By (6.5) in Zhang (1988), there exists a constant $K^* < \infty$, such that

$$\sum_{n=n_0}^{\infty} n^{p-1} P_\nu \{T' > \tau' + n\}$$

$$\leq \sum_{n=n_0}^{\infty} n^{p-1} P_\nu \{S_\tau + n - S_{\tau+n} \geq \gamma' n\} + \sum_{n=n_0}^{\infty} n^{p-1} P_\nu \left\{ \max_{n_1 \leq j \leq n^*} U_j^{-1} > \gamma' n \right\}$$

$$+ \sum_{n=n_0}^{\infty} n^{p-1} P_\nu \{K^* (\tau - b)^2 / b > \gamma' n\}$$

$$+ \sum_{n=n_0}^{\infty} n^{p-1} P_\nu \{K^* |\tau - b| n_1^{-1/2} > \gamma' n\} + o(1).$$

It follows from Lemma 1(i), (ii), (3.21) and the condition $E_\nu \{|\xi_1|^{2pp'} V(X_0)\} < \infty$ for some $p' > 1$ that

$$\sum_{n=n_0}^{\infty} n^{p-1} P_\nu \{T > \tau' + n\} \to 0 \quad \text{as } \min(n_0, b) \to \infty.$$

This proves the uniform integrability of $\{(T - \tau)^{+p}\}$ since the uniform integrability of $\{T^p; b \leq b^*, \lambda \in \Lambda\}$ for any given $b^*$ is implied by (6.2) in Zhang (1988). □

PROOF OF THEOREM 2. (i) For ease of notation, let $T = T_\lambda$, $n_1 = [b - \gamma b^\alpha]$, $n_2 = [b + \gamma b^\alpha]$, $T' = \max(n_1, \min(T, n_2))$ and $\tau' = \max(n_1, \min(\tau, n_2))$. Though $n_2$ is different from $n^*$ in Lemma 3, $(T - T')^+ \leq (T - \tau)^+ + (\tau - n_2)^+$, and by Lemma 4 we have

$$\lim_{b\to\infty} E_\nu |T - T'|^p = \lim_{b\to\infty} E_\nu |\tau - \tau'|^p = 0. \tag{6.3}$$

Clearly,

$$P_\nu \{\tau' > T' + n\}$$
$$\leq P_\nu \{L \geq n_1\} + P_\nu \{T \leq n_1\} + P_\nu \{\tau \geq n_2\} \tag{6.4}$$
$$+ P_\nu \{L < n_1 < T \leq T + n < \tau < n_2\}.$$

By (3.19) and Lemma 2,

$$\sum_{n=n_0}^{\infty} n^{p-1} P_\nu\{T' > \tau' + n\} = \sum_{n=n_0}^{n_2-n_1} n^{p-1} P_\nu\{T' > \tau' + n\}$$
$$\leq \sum_{n=n_0}^{n_2-n_1} n^{p-1} P_\nu\{L < n_1 < T \leq T+n < \tau < n_2\} + o(1). \quad (6.5)$$

On the event $\{L < n_1 < T \leq T+n < \tau < n_2\}$,

$$S_{T+n} + f(b) \leq b + d(T+n-b)$$
$$\leq \mu^* n + (b + d(T-b) - A(T;\lambda)) + S_T + U_T + f(T), \quad (6.6)$$

and by an argument similar to (6.10) of Zhang (1988), there exists a finite constant $K^*$ that does not depend on $\gamma$, and $\delta^* = \gamma K^* b^{\alpha-1} + K^* b^{-1/2}$, and by (6.5),

$$(6.7) \quad S_{T+n} - S_T \leq \mu^* n + U_T + \delta^* |\tau' - T'| + K^* (\tau-b)^2/b + K^* |\tau-b| b^{-1/2}.$$

Therefore, it follows from (6.5) and (6.7) that for $\gamma' = (1-\mu^*)/5$,

$$\sum_{n=n_0}^{\infty} n^{p-1} P_\nu\{\tau' > T' + n\}$$
$$\leq \sum_{n=n_0}^{\infty} n^{p-1} P_\nu\{S_T + n - S_{T+n} \geq \gamma' n\} + \sum_{n=n_0}^{\infty} n^{p-1} P_\nu\left\{\max_{n_1 \leq j \leq n_2} U_j \geq \gamma' n\right\}$$
$$(6.8)$$
$$+ \sum_{n=n_0}^{\infty} n^{p-1} P_\nu\{K^*(\tau-b)^2/b + K^*|\tau-b|n^{-1/2} \geq \gamma' n\}$$
$$+ \sum_{n=n_0}^{\infty} n^{p-1} P_\nu\{\delta^*(\tau'-T') \geq \theta' n\} + o(1) \quad \text{as } \min(n_0, b) \to \infty.$$

Since $\gamma$ is arbitrary, we can choose $\gamma$ small enough such that $(4\delta^*/\gamma')^p \leq 2$. Hence, it follows from (6.8), Lemma 1(i), (ii), (3.21) and Lemma 4 that as $\min(n_0, b) \to \infty$,

$$(6.9) \quad \sum_{n=n_0}^{\infty} n^{p-1} P_\nu\{\tau' - T' > n\} \leq \sum_{n=n_0}^{\infty} n^{p-1} P_\nu\{\delta^*(\tau'-T') \geq \gamma n\} + o(1).$$

And by (6.3) and Lemma 4, $\sum_{n=n_0}^{\infty} n^{p-1} P_\nu\{\tau' - T' > n\} = o(1)$, and

$$\{|T-\tau|^p; \lambda \in \Lambda\} \text{ is uniformly integrable.}$$

(ii) For the case where $\partial^2 A/\partial t^2 = 0$, the term $(\tau-b)^2/b$ disappears throughout the proof of (i). □

PROOF OF THEOREM 3. By definition, $S_T = R + A(T;\lambda) - \eta_T$ and $S_\tau = R(c_\lambda, d) + \mu b + d(\tau - b) - f(b)$. It follows that

$$S_T - S_\tau - d(T - \tau) \tag{6.10}$$
$$= R - R(c_\lambda, d) + [A(T;\lambda) - \mu b - d(T - b)] - [\eta_T - f(b)].$$

Recall $\Delta$ defined in (3.28). Let $d_1 = \xi_1 - \mu + \Delta(x_1; g) - \Delta(x_0; g)$, $\tilde{g} = E(d_1^2|X_0, X_1)$, $\sigma^2 = E_\pi d_1^2$ and $\Delta_2(x; g) = \Delta(x; \tilde{g})$. Since

$$S_T - \mu T - S_\tau + \mu \tau$$
$$= [(S_{\max(T,\tau)} - \mu \max(T, \tau)) - (S_\tau - \mu \tau)]$$
$$- [(S_\tau - \mu \tau) - (S_{\min(T,\tau)} - \mu \min(T, \tau))],$$

and A1–A3 imply that the conditions in Fuh and Zhang (2000) hold as shown in Lemma 1, it follows from Markov Wald's equation for Markov chains for second moments in Corollary 1 of Fuh and Zhang (2000) that

$$E_\nu(S_T - S_\tau - \mu(T - \tau))^2$$
$$= \sigma^2[E_\nu(\max(T, \tau) - \tau) + E_\nu(\tau - \min(T, \tau))]$$
$$- 2E_\nu\{(S_T - S_\tau - \mu(T - \tau))\Delta(X_{|T-\tau|})\}$$
$$+ E_\nu\{\Delta_2(X_{|T-\tau|}) - \Delta_2(X_0)\}$$
$$= (\sigma^2 - 2\mu)E_\nu|T - \tau| + O(1).$$

Therefore, by Theorem 2, $S_T - S_\tau - d(T - \tau)$ is uniformly integrable. By Anscombe's type central limit theorem for a Harris recurrent Markov random walk [Theorem 1 of Malinovskii (1986)], $\tau - b/\sqrt{b} \to (\mu - d_1^*)\sigma N(0, 1)$ in distribution. It follows that, as $b \to \infty$,

$$\frac{T - b}{\sqrt{b}} \to \sigma^* N(0, 1) \text{ in distribution by } (3.27), \tag{6.11}$$

$$\lim_{b \to \infty} P_x\{R > r\} = \lim P_x\{R(c, d) > r\} = G(r, d^*) \quad \text{by Theorem 1,} \tag{6.12}$$

$$A(T;\lambda) - \mu b - d(T - b) \to d_2^*(\sigma^* N(0, 1))^2/2 \tag{6.13}$$
$$\text{in distribution by } (3.34),$$

$$\eta_T - f(b) \to U \text{ in distribution by } (3.14), (3.31), (3.32), \tag{6.14}$$

where $c = c_\lambda = (\mu - d_\lambda)b_\lambda - f(b_\lambda)$ and $\sigma^* = (\mu - d_1^*)^{-1}\sigma$.

By an argument similar to Theorem 3(i) of Fuh and Lai (1998), $R(c_\lambda, d_\lambda)$ is uniformly integrable. Hence,

$$E_\nu(S_T - S_\tau - d(T - \tau)) = (\mu - d_1^*)^{-2}\sigma^2 d_2^*/2 - E_\nu U + o(1),$$
$$E_\nu T = E_\nu \tau + (\mu - d_1^*)^{-3}\sigma^2 d_2^*/2 - (\mu - d_1^*)^{-1} E_\nu U + o(1),$$

and

$$E_\nu \tau = (\mu - d)^{-1} c + (\mu - d^*)^{-1} \left( r(d^*) - \int \Delta(x) \, d(\pi_+^{d^*}(x) - \nu(x)) \right) + o(1)$$

$$= b_\lambda - (\mu - d)^{-1} f(b_\lambda)$$

$$+ (\mu - d_1^*)^{-1} \left( r(d_1^*) - \int \Delta(x) \, d(\pi_+^{d^*}(x) - \nu(x)) \right) + o(1).$$

This completes the proof. □

**7. Proofs of Theorems 4 and 5.** To establish asymptotic optimality of the SRP rule and derive the second-order asymptotic approximation for the average run length, we need to apply nonlinear Markov renewal theory developed in Section 3. Note that the Markov chain $\{W_n^\theta, n \geq 0\}$ on $\mathcal{X} := D' \times P(R^d) \times P(R^d)$ is induced by the products of random matrices $\{M_n, n \geq 0\}$. A positivity hypothesis on the matrices in the support of the Markov chain leads to contraction properties, on which basis the spectral theory is developed in Fuh (2003). Another natural hypothesis is that the transition probability possesses a density. This leads to a classical situation in the context of the so-called "Doeblin condition" for Markov chains. It also leads to precise results of the limiting theory and has been used to prove a nonlinear renewal theory in Section 3. We summarize the behavior of $\{W_n^\theta, n \geq 0\}$ in the following proposition. Note that in the case of i.i.d. iterated random functions satisfying Lipschitz conditions, similar results can be found in Theorems 2.1, 2.2 and Corollary 2.3 of Alsmeyer (2003). Here we generalize it to Markovian products of random matrices.

PROPOSITION 2. *Consider a given hidden Markov chain as in* (2.1) *and* (2.2) *satisfying* C1 *and* C2, *and let* $\theta = (\theta_0, \theta_1) \in \Theta \times \Theta$ *be the parameters. Then the induced Markov chain* $\{W_n^\theta, n \geq 0\}$ *defined in* (2.13) *is an aperiodic, $\mu$-irreducible and Harris recurrent Markov chain. Moreover, it is also a $V$-uniformly ergodic Markov chain for some $V$ on $\mathcal{X}$. And we have* $\sup_w \{\mathbb{E}_w(V(W_1))/V(w)\} < \infty$, *and there exist* $a, C > 0$ *such that* $\mathbb{E}_w(\exp\{a\chi(M_1)\}) \leq C$ *for all* $w = (y, \bar{u}, \bar{u}) \in \mathcal{X}$.

PROOF. For simplicity of notation, we delete $\theta$ in $\{W_n^\theta, n \geq 0\}$ in the proof. First, we prove that $\{W_n, n \geq 0\}$ is Harris recurrent. Note that the transition probability kernel of the Markov chain $\{(X_n, \xi_n), n \geq 0\}$ defined in (2.2) has probability density function, and the random matrices defined in (2.9) and (2.10) also have probability density with respect to $\mu$. Therefore, there exists a measurable function $g: \mathcal{X} \times \mathcal{X} \to [0, \infty)$ such that

(7.1) $$\mathbb{P}(w, dw') = g(w, w') \, d\mu(w'),$$

where $\int_{\mathcal{X}} g(w, w') \, d\mu(w') > 0$ for all $w \in \mathcal{X}$. For an arbitrary stopping time $\tau = h(W_n)$ for $W_n$, let $\mathbb{P}^\tau(w, \cdot) := \mathbb{P}_w(W_\tau \in \cdot)$ for $w \in \mathcal{X}$. For $A \in \mathcal{B}(\mathcal{D}')$ and $B \in \mathcal{B}(P(R^d) \times P(R^d))$, define

$$\mu^\tau(A \times B) := \int_{\mathcal{X}} \mathbb{P}\{W_\tau(w') \in A \times B\} \, d\mu(w').$$

Then

$$\mathbb{P}^{\tau+1}(w, A \times B) = \int_{\mathcal{X}} \mathbb{P}^\tau(w', A \times B) g(w, w') \, d\mu(w')$$
$$= \int_{\mathcal{X}} \mathbb{P}\{W_\tau(w') \in A \times B\} g(w, w') \, d\mu(w')$$

for all $A \in \mathcal{B}(\mathcal{D}')$ and $B \in \mathcal{B}(P(R^d) \times P(R^d))$. Therefore, given any $\mathbb{P}_m$-a.s. finite stopping time $\tau$ for $\{W_n, n \geq 0\}$, the family $(\mathbb{P}^{\tau+1}(w, \cdot))_{w \in \mathcal{X}}$ is nonsingular with respect to $\mu^\tau$.

We have thus particularly shown that, if $\mathbb{P}$ has a probability density with respect to $\mu$, then $\mathbb{P}^n$ has a probability density with respect to $\mu$ for all $n \geq 1$ (with, in general, different $\mu$). Let $g_\tau$ be such that

(7.2) $\qquad \mathbb{P}^{\tau+1}(w, dw') = g_\tau(w, w') \, d\mu^\tau(w'), \qquad w \in \mathcal{X},$

where $\int_{\mathcal{X}} g_\tau(w, w') \, d\mu^\tau(w') > 0$ for all $w \in \mathcal{X}$. It is easy to check that all $\mu$ and $\mu^\tau$ are absolutely continuous with respect to $m$.

Next, under condition C1, for each $m$-positive $A \times B$ let

$$\Gamma_0(A \times B) := \{w \in \mathcal{X} : \mathbb{P}_w\{W_n \in A \times B \text{ i.o.}\} = 1\}$$

satisfy $m(\Gamma_0(A \times B)) = 1$ and, thus, also $\mathbb{P}(w, \Gamma_0(A \times B)) = 1$ for $m$-almost all $w \in \mathcal{X}$. Recursively, define

$$\Gamma_{n+1}(A \times B) := \{w \in \Gamma_n(A \times B) : \mathbb{P}(w, \Gamma_n(A \times B)) = 1\}$$

for $n \geq 0$. Then $m(\Gamma_n(A \times B)) = 1$ for all $n \geq 0$ and $\Gamma_n(A \times B) \downarrow \Gamma_\infty(A \times B) := \bigcap_{k \geq 0} \Gamma_k(A \times B)$, as $n \to \infty$, giving $m(\Gamma_\infty(A \times B)) = 1$. Furthermore, $\Gamma_\infty(A \times B)$ is absorbing because, by construction, $\mathbb{P}(w, \Gamma_n(A \times B)) = 1$ for all $w \in \Gamma_\infty(A \times B)$ and $n \geq 0$, and, thus, $\mathbb{P}(w, \Gamma_\infty(A \times B)) = \lim_{n \to \infty} \mathbb{P}(w, \Gamma_n(A \times B)) = 1$ for all $w \in \Gamma_\infty(A \times B)$.

In particular, put $\tau = 1$. Denote $B^c$ as the complement of $B$. Since $m(\Gamma_\infty(\mathcal{X})^c) = 0$, also $\mu(\Gamma_\infty(\mathcal{X})^c) = 0$. It is now obvious from the previous considerations that we can choose $\delta > 0$ sufficiently small such that

$$\int_{\Gamma_\infty(\mathcal{X})} \int_{\mathcal{X}} \int_{\Gamma_\infty(\mathcal{X})} \mathbb{1}_{\{g \geq \delta\}}(w_1, w_2) \mathbb{1}_{\{g \geq \delta\}}(w_2, w_3) \, d\mu(w_3) \, d\mu(w_2) \, dm(w_1) > 0.$$

Hence, by Lemma 4.3 of Niemi and Nummelin (1986), there exist an $m$-

positive set $\Gamma_1 \subset \Gamma_\infty(\mathcal{X})$ and a $\mu$-positive set $\Gamma_2 \subset \Gamma_\infty(\mathcal{X})$ such that

$$\alpha := \inf_{w_1 \in \Gamma_1, w_3 \in \Gamma_2} \mu^2 \{w_2 \in \mathcal{X} : g(w_1, w_2) \geq \delta, \ g(w_2, w_3) \geq \delta\} > 0.$$

A combination of the above result with (7.1) and (7.2) implies

$$\begin{aligned}
\mathbb{P}^3(w_1, A \times B) &= \int_\mathcal{X} \mathbb{P}(w_2, A \times B) \mathbb{P}^2(w_1, dw_2) \\
&\geq \int_\mathcal{X} g_2(w_1, w_2) \int_{A \times B \cap \Gamma_2} g(w_2, w_3) \, d\mu(w_3) \, d\mu(w_2) \\
&\geq \alpha \delta^2 \mu(A \times B \cap \Gamma_2)
\end{aligned} \tag{7.3}$$

for all $w_1 \in \Gamma_1$ and $A \times B \in \mathcal{B}(\mathcal{X})$. By defining $\mathbb{H} := \Gamma_\infty(\mathcal{X}_1)$, we obtain an absorbing set such that $\Gamma_1$ is a regeneration set for $\{W_n, n \geq 0\}$ restricted to $\mathbb{H}$, that is, $\Gamma_1$ is recurrent and satisfies a minorization condition, namely (7.3). This proves the Harris recurrence of $\{W_n, n \geq 0\}$ on $\mathbb{H}$. By the previous construction, it is easy to see that $\mathbb{H} = \mathcal{X}$. Since $\{W_n, n \geq 0\}$ possesses a stationary distribution, it is clearly positive Harris recurrent.

Next, we give the proof of aperiodicity. If $\{W_n, n \geq 0\}$ were $q$-periodic with cyclic classes $\Gamma_1, \ldots, \Gamma_q$, say, then the $q$-skeleton $(W_{nq})_{n \geq 0}$ would have stationary distributions $\frac{m(\cdot \cap \Gamma_k)}{m(\Gamma_k)}$ for $k = 1, \ldots, q$. On the other hand, $Y_n$ is aperiodic by definition and $\overline{T_{nq}(\theta)u}$ is also a product of random matrices satisfying condition C1 and thus possesses only one stationary distribution. Consequently, $q = 1$ and $\{W_n, n \geq 0\}$ is aperiodic.

Note that we have $\mathbb{P}_w\{W_n \in A \times B \text{ i.o.}\} = 1$ for all $w \in \mathcal{X}$ and all $m$-positive open $A \times B \in \mathcal{B}(\mathcal{X})$. Denote $\rho(B)$ as the first return time to $B$ for $W_n$. Hence, $m(\text{int}(\mathcal{X})) > 0$ ensures $\mathbb{P}_w(\rho(\mathcal{X}) < \infty) = 1$ for all $w \in \mathcal{X}$, which easily yields the $\mu$-irreducibility of $\{W_n, n \geq 0\}$.

Under conditions C1 and C2, the property of $V$-uniform ergodicity is taken from Lemma 4 of Fuh (2003). The finiteness A1–A3 of the moments comes from C2 and a simple calculation. The details are omitted. □

To prove the main results in Section 4, our first aim is to find a sequence of $p$'s that converge to 0, for which the stopping times $N_{B,p}$ converge to an appropriate stopping time. Furthermore, for technical reasons, we want all the stopping times in the sequence to be bounded by some stopping time with finite expectation.

LEMMA 5. *Consider the problem $B(\beta = 0, p, c, \tilde{w})$ described in Section 4. Then the following hold:*

(i) *There exist a constant $D_c$ and some $0 < q_0 < 1$, such that for all $q_0 \leq q \leq 1$ and for all threshold functions $B(\cdot)$ with the property that $B(w) \geq D_c$, for each $w \in \mathcal{X}$, we have*

$$(7.4) \qquad \mathbb{E}_w(N_{B,p} - \omega | N_{B,p} \geq \omega) > 2c^{-1},$$

where $N_{B,p}$ is defined as in (4.6).

(ii) Define $N_{D,p,w} = \inf\{n : LR_{n,p} \geq D, W_n = w\}$. Then for each $0 \leq p \leq 1 - q_0$, there exists $w = w(p)$ such that with probability 1,

$$(7.5) \qquad N_{B,p} \leq N_{D,p,w} \leq N_{D,1,w}.$$

Furthermore, there exist a state $w_1 \in \mathcal{X}$ and a subsequence of $p$'s, such that (7.5) is true with $w(p) = w_1$ for all the $p$'s in the subsequence.

(iii) Let $B_p(\cdot)$ be the threshold function of the problem $B(\beta, p, c, \tilde{w})$, and assume that $B_p(\cdot) \to_{p \to 0} B_0(\cdot)$ for some function $B_0(\cdot)$. Assume further that the convergence is along the subsequence of $p$'s from (ii). Then $B_0(w_1) \leq D_c$.

(iv) Denote $D_c^* = \inf\{D_c | D_c \text{ as in (i)}\}$. Then $D_c^*$ is nonincreasing in $c$ and $D_c^* \to B_0(w_1)/q_0$ as $c \to \infty$.

PROOF. Since the induced Markov chain $\{W_n^\theta, n \geq 0\}$ is Harris recurrent via Proposition 2, we may assume, without loss of generality, that there exists a recurrent state $w_0$ of the Markov chain governed by $\mathbb{P}^{\theta_0}$. By making use of the regeneration scheme for Harris recurrent Markov chains, the proof of Lemma 5 is similar to that of Lemma 1 in Yakir (1994). The details are omitted. □

REMARK 3. Notice first that the constant $D_c$ does not depend on the initial state $(p, \tilde{w})$. Lemma 5 remains true when $c = c(p)$ is allowed to vary with $p$, as long as $\liminf_{p \to 0} c(p) > 0$. In particular, it is correct if $c(p)$ converges to some positive $c$.

Let $(1 - \rho(N, \nu))/p$ be the normalized risk of a stopping time $N$. Using the results of Lemma 5, we can show that for $p \to 0$, the (normalized) risk of a converging sequence of stopping times goes to a limit. Consider the Bayesian problem $B(\beta = 0, p, c, \tilde{w})$, and let $N$ be a stopping time. A similar argument to that of Lemma 9 in Pollak (1985) implies that as $p \to 0$,

$$(7.6) \qquad \frac{P_\pi(N > \omega)}{p} \to E_\infty N.$$

LEMMA 6. Let $B_p(\cdot)$ be defined as in (4.6) and let $e_c(\cdot) = \liminf_{p \to 0} B_p(\cdot)/p$. Then with probability one, $\liminf_{c \to 0} e_c(\cdot) = \infty$.

PROOF. Let $N_{B,p}$ be defined as in (4.6). Suppose for all $w \in \mathcal{X}$, $\liminf_{c \to 0} e_c(w) = e_\infty(w) < \infty$. Then for almost all $w \in \mathcal{X}$, $B_{p_i}(w)/(p_i(1 - B_{p_i}(w))) < 1 + e_\infty(w)$ for some subsequence $c_i \to 0, p_i \to 0$ as $i \to \infty$. Since

$$E_\pi(\text{Loss using } N_{B,p})$$
$$= P_\pi(N_{B,p} < \omega) + cP_\pi(N_{B,p} \geq \omega)E_\pi(N_{B,p} - \omega | N_{B,p} \geq \omega),$$

it follows from (7.6) that

$$\frac{1 - E_\pi(\text{Loss using } N_{B,p})}{p_i}$$
$$= \frac{P_\pi(N_{B,p} \geq \omega)}{p_i}(1 - c_i E_\pi(N_{B,p} - \omega | N_{B,p} \geq \omega))$$
$$\leq \frac{P_\pi(N_{B,p} \geq \omega)}{p_i} \leq \frac{P_\pi(N_{1+e_\infty} \geq \omega)}{p_i} \leq 1 + E_\infty N_{1+e_\infty}$$

for large enough $i$. Clearly, $E_\infty N_{1+e_\infty} < \infty$. Hence, one can do better by using a CUSUM rule in the hidden Markov model with large enough upper boundary [Fuh (2003)], and this contradicts the fact that $N_{B,p}$ is a Bayes rule. □

LEMMA 7. *Let $N_b^\psi$ be defined as in (4.14), and $N_{B,p}$ be defined as in (4.6). Assume the boundary $B(w)$ defined in (4.6) is chosen as $Bg(w)$ for a measurable function $g$ with $\int_\mathcal{X} g(w)\,dm(w) < \infty$. Then $\mathbb{E}_1(N_{B,p}|W_0 = \tilde{w}) = \mathbb{E}_1(N_b^\psi|W_0 = \tilde{w}) + o(1)$ as $p \to 0$ and $B \to \infty$.*

The proof of Lemma 7 is given in Section 8.

PROOF OF THEOREM 4. By using Proposition 2 and Lemmas 5–7, the proof of Theorem 4 is similar to that of Theorem 1 in Pollak (1985). The details are omitted. □

LEMMA 8. *For each $B$, let $T_B$ be defined as in (4.11) and (4.12). Then $T_B F_n = F_{n+1}$, and, hence, there associates a set of invariant measures $\Phi_B$ such that $T_B \phi = \phi$ for all $\phi \in \Phi_B$.*

PROOF. Since

$$\mathbb{P}_\infty(R_{n,p} \in ds, W_{n+1} \in dw' | N_{q,b} > n+1, W_n = w)$$
$$= \mathbb{P}_\infty(R_{n,p} \in ds, W_{n+1} \in dw' | N_{q,b} > n, N_{q,b} > n+1, W_n = w)$$
$$= \mathbb{P}_\infty(N_{q,b} > n+1, W_{n+1} \in dw' | R_{n,p} = t, N_{q,b} > n, W_n = w)$$
$$\times \left[ \int_{w,w' \in \mathcal{X}} \int_0^{B(w')} \mathbb{P}_\infty(N_{q,b} > n+1, W_{n+1} \in dw' | R_{n,p} = t, \right.$$
$$\left. N_{q,b} > n, W_n = w) \right]^{-1}$$
$$\times \frac{\mathbb{P}_\infty(R_{n,p} \in ds | N_{q,b} > n, W_n = w)\mathbb{P}_\infty(N_{q,b} > n | W_n = w)\mathbb{P}(w, dw')}{\mathbb{P}_\infty(R_{n,p} \in ds | N_{q,b} > n, W_n = w)\mathbb{P}_\infty(N_{q,b} > n | W_n = w)\mathbb{P}(w, dw')}$$

$$= \frac{\zeta(s,w,w')\mathbb{P}(w,dw')\,dF_n(s,w)}{\int_{w,w'\in\mathcal{X}}\int_0^{B(w')}\zeta(s,w,w')\mathbb{P}(w,dw')\,dF_n(s,w)},$$

it follows that

$$F_{n+1}(s,w) = \frac{\int_{w'\in\mathcal{X}}\int_0^{B(w')}\rho(t,s,w,w')\zeta(t,w,w')\mathbb{P}(w,dw')\,dF_n(t,w)}{Q(F_n)}$$
$$= T_B F_n(s,w).$$

The existence of the fixed point follows by the same argument as that of Lemma 11 in Pollak (1985). □

For a given $\gamma$, let $\mathcal{N}_\gamma$ be the set of all detection policies $N_b^\psi$, defined as in (4.14), for which $\mathbb{E}_\infty N_b^\psi = \gamma$. In the next lemma it is shown that $\mathcal{N}_\gamma$ is not empty. Furthermore, this set contains a stopping rule that is a limit of Bayes stopping rules.

LEMMA 9. *There exist a sequence of p's that converges to 0, a sequence of randomized Bayes problems $B(\beta=0,p,c(p),\psi_p)$ with the appropriate Bayes rules $N_{q,b}^{\psi_p}$, defined as the detection policy $N_b^\psi$ in (4.14) with $R_n$ replaced by $R_{q,n}$, and a constant $0 < c < \infty$ such that*

(i) $c(p) \to_{p\to 0} c$,
(ii) $\mathbb{E}_\infty N_b^\psi = \gamma$,
(iii) $\rho^{\psi_p}(N(\beta=0,p,c(p),\psi_p)) \to_{p\to 0} (\mu(\psi)+\gamma)(1-c\mathbb{E}_1 N_b^\psi)$,

*where $\mu(\psi) = \int_{w\in\mathcal{X}}\int_0^\infty r\,d\psi(r,w)$, and $\rho^{\psi_p}(N(\beta=0,p,c(p),\psi_p))$ is the normalized Bayes risk.*

PROOF. The proof is similar to that of Lemma 2 in Yakir (1994) and is omitted. □

By Lemma 7, the difference of the expected values for the stopping rule with constant boundary and the stopping rule with curved boundary is $o(1)$. Therefore, we only need to consider the stopping rule with a constant boundary in the following lemmas. Note that in this case, $\psi(s,w) = \psi(s)$.

LEMMA 10. *Let*

$$m(k) = \begin{cases} \dfrac{\int_0^B s\,d\psi(s)}{\int_0^B s\,d\psi(s) + \mathbb{E}_\infty N_b^\psi}, & k = 0, \\ \dfrac{\mathbb{P}_\infty(N_b^\psi \geq k)}{\int_0^B s\,d\psi(s) + \mathbb{E}_\infty N_b^\psi}, & k = 1,2,\ldots. \end{cases}$$

If one uses $N_{q_i,b}^{\psi_{p_i}}$ in the problem $B(\beta=0,p_i,c(p_i),\psi_{p_i})$, then $\mathbb{P}(\omega=k|N_{q_i,b}^{\psi_{p_i}}\geq\omega)\to m(k)$ as $i\to\infty$. Also, if one uses $N_b^{\psi}$ instead of $N_{q_i,b}^{\psi_{p_i}}$ in problem $B(\beta=0,p_i,c(p_i),\psi_{p_i})$, then $\mathbb{P}(\omega=k|N_b^{\psi}\geq\omega)\to m(k)$ as $i\to\infty$.

PROOF. Note that $\mathbb{P}(\omega=0)/p_i\to\int_0^B s\,d\psi(s)$ as $i\to\infty$. If one uses $N_{q_i,b}^{\psi_{p_i}}$ in problem $B(\beta=0,p_i,c(p_i),\psi_{p_i})$, then by using an argument similar to the proof of Lemmas 9 and 12 of Pollak (1985), we have

$$\text{(7.7)} \qquad \frac{\mathbb{P}(N_{q_i,b}^{\psi_{p_i}}\geq\omega)}{p_i} \xrightarrow[i\to\infty]{} \int_0^B s\,d\psi(s)+\mathbb{E}_\infty N_b^\psi.$$

Therefore,

$$\mathbb{P}(\omega=k|N_{q_i,b}^{\psi_{p_i}}\geq\omega)=\frac{\mathbb{P}(N_{q_i,b}^{\psi_{p_i}}\geq k|\omega=k)\mathbb{P}(\omega=k)}{\mathbb{P}(N_{q_i,b}^{\psi_{p_i}}\geq\omega)}\xrightarrow[i\to\infty]{} m(k).$$

A similar argument applies when using $N_b^\psi$ instead of $N_{q_i,b}^{\psi_{p_i}}$. □

By using the same argument as in Lemma 13 of Pollak (1985), we also have

$$\lim_{i\to\infty}\frac{1-\{\text{Expected loss using }N_b^\psi\text{ for }B(\beta=0,p_i,c(p_i),\psi_{p_i})\}}{p_i}$$
$$=\lim_{i\to\infty}\frac{1-\{\text{Expected loss using }N_{q_i,b}^{\psi_{p_i}}\text{ for }B(\beta=0,p_i,c(p_i),\psi_{p_i})\}}{p_i}$$
$$\text{(7.8)} \qquad \times\left[\int_0^B s\,d\psi(s)+\mathbb{E}_\infty N_b^\psi\right]$$
$$\times\left\{1-c^*\left[\mathbb{E}_1(N_b^\psi-1)\frac{\mathbb{E}_\infty N_b^\psi}{\int_0^B s\,d\psi(s)+\mathbb{E}_\infty N_b^\psi}\right.\right.$$
$$\left.\left.+\frac{\int_0^B s\,d\psi(s)}{\int_0^B s\,d\psi(s)+\mathbb{E}_\infty N_b^\psi}\lim_{i\to\infty}\mathbb{E}_\infty(N_b^\psi|\omega=0)\right]\right\}.$$

The following lemma generalizes Theorem 5 of Kesten (1973) from products of i.i.d. random matrices to products of Markov random matrices.

LEMMA 11. *Let $M_n$, $n\geq 0$ be the random matrices defined in (2.8) and (2.9). Assume that $\mathbb{E}_m\log|M_1|<0$, but that for some $k_1>0$, $\mathbb{E}_m|M_1|^{k_1}=1$, $\mathbb{E}_m|M_1|^{k_1}\times\log^+|M_1|<\infty$. Assume, in addition, $\log|M_1|$ does not have a lattice distribution. Then the series $R=\sum_{k=0}^\infty M_k\cdots M_1 M_0$ converges with probability 1,*

*and*

$$\lim_{t\to\infty} t^{k_1} \mathbb{P}_m(R > t) \quad and \quad \lim_{t\to\infty} t^{k_1} \mathbb{P}_m(R < -t)$$

*exist and are finite.*

PROOF. By making use of the result that for all $w \in \mathcal{X}$,

$$B\mathbb{P}_w\Big\{\max_{n\geq 1} \|M_n \cdots M_1 M_0 \pi\| > B\Big\} \to C \quad \text{as } B \to \infty,$$

for some constant $C$, developed in Section 3.2 of Fuh (2004) for products of Markov random matrices, the proof of the remainder part is similar to that of Theorem 4 and Theorem 5 in Kesten (1973). The details are omitted. $\square$

LEMMA 12. $\int_0^B s\, d\psi(s)/[\int_0^B s\, d\psi(s) + \mathbb{E}_\infty N_b^\psi] = O((\log B)/B)$, *where* $O((\log B)/B)/((\log B)/B)$ *remains bounded as* $B \to \infty$.

PROOF. Denote by $\mathcal{F}_n$ the $\sigma$-algebra generated by $\{\xi_0, \xi_1, \ldots, \xi_n\}$. Since $R_{n+1}^* = \beta(W_{n-1}^\theta, W_n^\theta)(1 + R_n^*)$, it follows that $\mathbb{E}_\infty(R_{n+1}^*|\mathcal{F}_n) = 1 + R_n^*$, and, therefore, $R_n^* - n$ is a $\mathbb{P}_\infty$-martingale with expectation $\mathbb{E}_\infty(R_n^* - n) = \mathbb{E}_\infty R_0^* = \int_0^B s\, d\psi(s)$.

By using the optional sampling theorem, we have that $\int_0^\infty s\, d\psi(s) = \mathbb{E}_\infty R_{N_b^\psi}^* - \mathbb{E}_\infty N_b^\psi$. Therefore, $\int_0^\infty s\, d\psi(s) + \mathbb{E}_\infty N_b^\psi = \mathbb{E}_\infty R_{N_b}^* \geq B$.

It is easy to see that for all $n$,

(7.9)
$$\psi(s) = \mathbb{P}_\infty(R_n^* \leq s | N_b^\psi > n)$$
$$\geq \mathbb{P}_\infty(R_n^* \leq s) \underset{n\to\infty}{\to} \lim_{n\to\infty} \mathbb{P}_\infty(R_n \leq s) = \mathbb{P}(R \leq s).$$

Note that the limit in the above equation (7.9) follows from $R_n^* - R_n = R_0^* \exp\{\sum_{i=1}^n \mathbb{S}_i\} \to 0$ a.s. $P_\infty$ as $n \to \infty$. Hence,

$$\int_0^B s\, d\psi(s) = \int_0^B \Big[1 - \int_0^s d\psi(t)\Big] ds \leq \int_0^B \mathbb{P}(R > s)\, ds.$$

Under the conditions of Lemma 10, this implies the conditions of Lemma 11 hold with $k_1 = 1$. And by Lemma 11, $s\mathbb{P}(R > s) \to 1$ as $s \to \infty$.

It follows that $\lim_{B\to\infty} \int_0^B s\, d\psi(s)/\log B \leq 1$, from which Lemma 12 follows.

$\square$

PROOF OF THEOREM 5. Let $N_b^\psi$ be a stopping time from the set $\mathcal{N}_\beta$ that minimizes $\mathbb{E}_k N$ among all stopping times $N$ from that set. The change point detection policy $N_b^\psi$ is a minimax policy in the sense of equations (2.3) and (2.4). Notice that a limit of Bayes rules minimizes $\mathbb{E}_k N$ among all stopping times in the set $\mathcal{N}_\beta$, hence the claim of the above theorem is not empty. By (4.15), $N_b^\psi$ is an equalizer rule and note that

$$N_b^\psi \leq N_b \leq \min\left\{n \Big| \max_{1 \leq k \leq n} \frac{p_k(\xi_0, \xi_1, \ldots, \xi_k; \theta_1)}{p_k(\xi_0, \xi_1, \ldots, \xi_k; \theta_0)} \geq B\right\},$$

which is the CUSUM stopping rule for hidden Markov models. By Theorem 7 of Fuh (2003), we have $\mathbb{E}_1 N_b^\psi = O(\log B)$. The rest of the proof is almost identical to the proof of Theorem 2 in Pollak (1985) and is thus omitted. $\square$

## 8. Proofs of Theorem 6 and Lemma 7.

PROOF OF THEOREM 6. Note that the probability $\mathbb{P}_1$ and expectation $\mathbb{E}_1$ in this section are taken under $W_0 = \tilde{w}$, and we omit it for simplicity. The proof of (5.10) rests on the nonlinear Markov renewal theory from Theorem 3 and Corollary 1. Indeed, by (5.3), the stopping time $N_b^\psi$ is based on the thresholding of the sum of the Markov random walk $\mathbb{S}_n$ and the nonlinear term $\eta_n$. Note that

$$\eta_n \underset{n \to \infty}{\to} \eta \quad \mathbb{P}_1\text{-a.s.} \quad \text{and} \quad \mathbb{E}_1 \eta_n \underset{n \to \infty}{\to} \mathbb{E}_1 \eta,$$

and $\eta_n$, $n \geq 1$ are slowly changing under $\mathbb{P}_1$. In order to apply Theorem 3 and Corollary 1, we have to check the validity of the following three conditions:

(8.1) $\quad \sum_{n=1}^{\infty} \mathbb{P}_1\{\eta_n \leq -\varepsilon n\} < \infty \quad$ for some $0 < \varepsilon < K(\mathbb{P}^{\theta_1}, \mathbb{P}^{\theta_0})$;

(8.2) $\quad \max_{0 \leq k \leq n} |\eta_{n+k}|, \quad n \geq 1, \text{ are } \mathbb{P}_1\text{-uniformly integrable};$

(8.3) $\quad \lim_{b \to \infty} b\, \mathbb{P}_1\left\{N_b^\psi \leq \frac{\varepsilon b}{K(\mathbb{P}^{\theta_1}, \mathbb{P}^{\theta_0})}\right\} = 0 \quad$ for some $0 < \varepsilon < 1$.

Condition (8.1) holds trivially because $\eta_n \geq 0$. Since $\eta_n$, $n = 1, 2, \ldots$, are nondecreasing, $\max_{0 \leq k \leq n} |\eta_{n+k}| = \eta_{2n}$ and to prove (8.2) it suffices to show that $\eta_n$, $n \geq 1$, are $\mathbb{P}_1$-uniformly integrable. Since $\eta_n \leq \eta$ and, by (5.9), $\mathbb{E}_1 \eta < \infty$, the desired uniform integrability follows. Therefore, condition (8.2) is satisfied.

We now turn to checking condition (8.3). By using $E_\pi \xi_1 > 0$, and $0 < K(\mathbb{P}^{\theta_1}, \mathbb{P}^{\theta_0}) < \infty$, we will prove that

(8.4) $\quad \mathbb{P}_1\left\{N_b^\psi < \frac{(1-\varepsilon)b}{K(\mathbb{P}^{\theta_1}, \mathbb{P}^{\theta_0})}\right\} \leq e^{-y_\varepsilon b} + \alpha_1(\varepsilon, b),$

where $y_\varepsilon > 0$ for all $\varepsilon > 0$ and

$$\alpha_1(\varepsilon, b) = \mathbb{P}_1\left\{\max_{1 \leq n < K_{\varepsilon,b}} \mathbb{S}_n \geq (1+\varepsilon)(1-\varepsilon)b\right\},$$

(8.5)

$$K_{\varepsilon,b} = (1-\varepsilon)b/K(\mathbb{P}^{\theta_1}, \mathbb{P}^{\theta_0}).$$

If (8.4) is correct, then the first term on the right-hand side of (8.4) is $o(1/b)$ as $b \to \infty$. All it remains to do is to show that $\alpha_1(\varepsilon, b)$ in (8.5) is $o(1/b)$.

To this end, by Proposition 2 we can apply Theorem 2 of Fuh and Zhang (2000) to have that for all $\varepsilon > 0$ and $r \geq 0$,

(8.6) $$\sum_{n=1}^{\infty} n^{r-1} \mathbb{P}_1\left\{\max_{1 \leq k \leq n}(\mathbb{S}_k - K(\mathbb{P}^{\theta_1}, \mathbb{P}^{\theta_0})k) \geq \varepsilon n\right\} < \infty,$$

whenever $\mathbb{E}_1|\mathbb{S}_1|^2 < \infty$ and $\mathbb{E}_1[(\mathbb{S}_1 - K(\mathbb{P}^{\theta_1}, \mathbb{P}^{\theta_0}))^+]^{r+1} < \infty$. Recall that under the conditions of Theorem 6, $\mathbb{E}_1|\mathbb{S}_1|^2 < \infty$, and hence, the sum on the left-hand side of the inequality (8.6) is finite for $r = 1$ and all $\varepsilon > 0$, which implies that the summand should be $o(1/n)$. Since

$$\alpha_1(\varepsilon, b) \leq \mathbb{P}_1\left\{\max_{n < K_{\varepsilon,b}}(\mathbb{S}_n - K(\mathbb{P}^{\theta_1}, \mathbb{P}^{\theta_0})n) \geq \varepsilon(1-\varepsilon)b\right\},$$

it follows that $\alpha_1(\varepsilon, b) = o(1/b)$.

Next, we need to prove (8.4). Denote $\mathbb{S}_n^k = \log LR_n^k$, and let $N = N_b^\psi$ for simplicity. For any $C > 0$, by using a change of measure argument, we have

$$\mathbb{P}_\infty\{N < (1-\varepsilon)bK(\mathbb{P}^{\theta_1}, \mathbb{P}^{\theta_0})^{-1}\} = \mathbb{E}_1\{\mathbb{1}_{\{N < (1-\varepsilon)bK(\mathbb{P}^{\theta_1}, \mathbb{P}^{\theta_0})^{-1}\}} e^{-\mathbb{S}_N^k}\}$$

$$\geq \mathbb{E}_1\{\mathbb{1}_{\{N < (1-\varepsilon)bK(\mathbb{P}^{\theta_1}, \mathbb{P}^{\theta_0})^{-1}, \mathbb{S}_N^k < C\}} e^{-\mathbb{S}_N^k}\}$$

$$\geq e^{-C} \mathbb{P}_1\left\{N < (1-\varepsilon)bK(\mathbb{P}^{\theta_1}, \mathbb{P}^{\theta_0})^{-1}, \max_{n < (1-\varepsilon)bK(\mathbb{P}^{\theta_1}, \mathbb{P}^{\theta_0})^{-1}} \mathbb{S}_n^k < C\right\}$$

$$\geq e^{-C}\left[\mathbb{P}_1\{N < (1-\varepsilon)bK(\mathbb{P}^{\theta_1}, \mathbb{P}^{\theta_0})^{-1}\}\right.$$

$$\left. - \mathbb{P}_1\left\{\max_{n < (1-\varepsilon)bK(\mathbb{P}^{\theta_1}, \mathbb{P}^{\theta_0})^{-1}} \mathbb{S}_n^k \geq C\right\}\right].$$

Choosing $C \leq (1+\varepsilon)(1-\varepsilon)b$, we then have

$$\mathbb{P}_1\left\{N < \frac{(1-\varepsilon)b}{K(\mathbb{P}^{\theta_1}, \mathbb{P}^{\theta_0})}\right\}$$

(8.7)

$$\leq e^C \mathbb{P}_\infty\{N < (1-\varepsilon)bK(\mathbb{P}^{\theta_1}, \mathbb{P}^{\theta_0})^{-1}\} + \alpha_1(\varepsilon, b).$$

Recall that $R_n^*$ is defined in (4.13). Note that under the condition $0 < K(\mathbb{P}^{\theta_1}, \mathbb{P}^{\theta_0}) < \infty$, we have

$$\mathbb{P}_\infty\{N < K_{\varepsilon,b}\} = \sum_{i=1}^{K_{\varepsilon,b}} \mathbb{P}_\infty\{R_i^* > B\} \leq \sum_{i=1}^{K_{\varepsilon,b}} \frac{i}{B} \leq \frac{(\log B)^2}{(K(\mathbb{P}^{\theta_1}, \mathbb{P}^{\theta_0}))^2 B}.$$

By choosing a suitable $C$, we have the first term of $(8.7) \leq e^{-y_\varepsilon b}$, for some $y_\varepsilon > 0$, and get the proof of (8.4).

Thus, all conditions of Theorem 3 are satisfied. The use of this theorem yields (5.10) for a large $b$. $\square$

PROOF OF LEMMA 7. Note that the parameter $\lambda$ defined in (3.2) is $B$ in this case. Since $B(w) = Bg(w)$ is independent of $t$ and $\int_\mathcal{X} g(w)\,dm(w) < \infty$, $d_B$ defined in (3.7) is 0. By (8.1)–(8.3) developed in the proof of Theorem 6, we can apply Corollary 1, to obtain that as $B \to \infty$,

$$\mathbb{E}_1(N_{B,p}|W_0 = \tilde{w})$$
$$(8.8) \quad = \frac{1}{E_\pi \xi_1(p)}\left(b + \frac{\mathbb{E}_{m_+}\mathbb{S}^2_{m_+}(p)}{2\mathbb{E}_{m_+}\mathbb{S}_{m_+}(p)} - \mathbb{E}_{m_+}\eta - \int \Delta_p(w)\,dm_+(w) + \Delta_p(\tilde{w})\right)$$
$$\quad + o(1).$$

Here the random variables in (8.8) are the corresponding terms of (2.15) divided by $q$. We also have $E_\pi\xi_1(p) \to E_\pi\xi_1$, $\mathbb{E}_{m_+}\mathbb{S}^2_{m_+}(p) \to \mathbb{E}_{m_+}\mathbb{S}^2_{m_+}$, $\mathbb{E}_{m_+}\mathbb{S}_{m_+}(p) \to \mathbb{E}_{m_+}\mathbb{S}_{m_+}$ and $\Delta_p(w) \to \Delta(w)$. By using Corollary 1 again, we have

$$\mathbb{E}_1(N_b^\psi|W_0 = \tilde{w})$$
$$= \frac{1}{E_\pi \xi_1}\left(b + \frac{\mathbb{E}_{m_+}\mathbb{S}^2_{m_+}}{2\mathbb{E}_{m_+}\mathbb{S}_{m_+}} - \mathbb{E}_{m_+}\eta - \int \Delta(w)\,dm_+(w) + \Delta(\tilde{w})\right) + o(1).$$

Hence, Lemma 7 is proved. $\square$

## REFERENCES

ALSMEYER, G. (1994). On the Markov renewal theorem. *Stochastic Process. Appl.* **50** 37–56. MR1262329
ALSMEYER, G. (2003). On the Harris recurrence of iterated random Lipschitz functions and related convergence rate results. *J. Theoret. Probab.* **16** 217–247. MR1956829
BANSAL, R. K. and PAPANTONI-KAZAKOS, P. (1986). An algorithm for detecting a change in a stochastic process. *IEEE Trans. Inform. Theory* **32** 227–235. MR838411
BASSEVILLE, M. and NIKIFOROV, I. V. (1993). *Detection of Abrupt Changes*: *Theory and Application*. Prentice Hall, Englewood Cliffs, NJ. MR1210954
BRAUN, J. V. and MÜLLER, H. G. (1998). Statistical methods for DNA sequence segmentation. *Statist. Sci.* **13** 142–162.
FUH, C.-D. (2003). SPRT and CUSUM in hidden Markov models. *Ann. Statist.* **31** 942–977. MR1994736

Fuh, C.-D. (2004). Uniform Markov renewal theory and ruin probabilities in Markov random walks. *Ann. Appl. Probab.* **14** 1202–1241. MR2071421

Fuh, C.-D. and Lai, T. L. (1998). Wald's equations, first passage times and moments of ladder variables in Markov random walks. *J. Appl. Probab.* **35** 566–580. MR1659504

Fuh, C.-D. and Lai, T. L. (2001). Asymptotic expansions in multidimensional Markov renewal theory and first passage times for Markov random walks. *Adv. in Appl. Probab.* **33** 652–673. MR1860094

Fuh, C.-D. and Zhang, C.-H. (2000). Poisson equation, maximal inequalities and $r$-quick convergence for Markov random walks. *Stochastic Process. Appl.* **87** 53–67. MR1751164

Götze, F. and Hipp, C. (1983). Asymptotic expansions for sums of weakly dependent random vectors. *Z. Wahrsch. Verw. Gebiete* **64** 211–239. MR714144

Kesten, H. (1973). Random difference equations and renewal theory for products of random matrices. *Acta Math.* **131** 207–248. MR440724

Kesten, H. (1974). Renewal theory for functionals of a Markov chain with general state space. *Ann. Probab.* **2** 355–386. MR365740

Lai, T. L. (1995). Sequential change point detection in quality control and dynamical systems (with discussion). *J. Roy. Statist. Soc. Ser. B* **57** 613–658. MR1354072

Lai, T. L. (1998). Information bounds and quick detection of parameter changes in stochastic systems. *IEEE Trans. Inform. Theory* **44** 2917–2929. MR1672051

Lai, T. L. (2001). Sequential analysis: Some classical problems and new challenges (with discussion). *Statist. Sinica* **11** 303–408. MR1844531

Lai, T. L. and Siegmund, D. (1977). A nonlinear renewal theory with applications to sequential analysis. I. *Ann. Statist.* **5** 946–954. MR445599

Lai, T. L. and Siegmund, D. (1979). A nonlinear renewal theory with applications to sequential analysis. II. *Ann. Statist.* **7** 60–76. MR515684

Lorden, G. (1971). Procedures for reacting to a change in distribution. *Ann. Math. Statist.* **41** 1897–1908. MR309251

Malinovskii, V. K. (1986). Asymptotic expansions in the central limit theorem for recurrent Markov renewal processes. *Theory Probab. Appl.* **31** 523–526. MR866885

Melfi, V. F. (1992). Nonlinear Markov renewal theory with statistical applications. *Ann. Probab.* **20** 753–771. MR1159572

Meyn, S. P. and Tweedie, R. L. (1993). *Markov Chains and Stochastic Stability.* Springer, New York. MR1287609

Moustakides, G. V. (1986). Optimal stopping times for detecting changes in distribution. *Ann. Statist.* **14** 1379–1387. MR868306

Niemi, S. and Nummelin, E. (1986). On non-singular renewal kernels with an application to a semigroup of transition kernels. *Stochastic Process. Appl.* **22** 177–202. MR860932

Pollak, M. (1985). Optimal detection of a change in distribution. *Ann. Statist.* **13** 206–227. MR773162

Pollak, M. (1987). Average run lengths of an optimal method of detecting a change in distribution. *Ann. Statist.* **15** 749–779. MR888438

Pollak, M. and Siegmund, D. (1975). Approximations to the expected sample size of certain sequential tests. *Ann. Statist.* **3** 1267–1282. MR403114

Ritov, Y. (1990). Decision theoretic optimality of the CUSUM procedure. *Ann. Statist.* **18** 1464–1469. MR1062720

Roberts, S. W. (1966). A comparison of some control chart procedures. *Technometrics* **8** 411–430. MR196887

Shiryayev, A. N. (1963). On optimum methods in quickest detection problems. *Theory Probab. Appl.* **8** 22–46.

Shiryayev, A. N. (1978). *Optimum Stopping Rules.* Springer, New York.

Siegmund, D. (1985). *Sequential Analysis. Tests and Confidence Intervals.* Springer, New York. MR799155

Woodroofe, M. (1976). A renewal theorem for curved boundaries and moments of first passage times. *Ann. Probab.* **4** 67–80. MR391294

Woodroofe, M. (1977). Second order approximations for sequential point and interval estimation. *Ann. Statist.* **5** 984–995. MR494735

Woodroofe, M. (1982). *Nonlinear Renewal Theory in Sequential Analysis.* SIAM, Philadelphia. MR660065

Yakir, B. (1994). Optimal detection of a change in distribution when the observations form a Markov chain with a finite state space. In *Change-Point Problems* (E. Carlstein, H. G. Müller and D. Siegmund, eds.) 346–358. IMS, Hayward, CA. MR1477935

Yakir, B. (1997). A note on optimal detection of a change in distribution. *Ann. Statist.* **25** 2117–2126. MR1474086

Zhang, C.-H. (1988). A nonlinear renewal theory. *Ann. Probab.* **16** 793–824. MR929079

Institute of Statistical Science
Academia Sinica
Taipei, 11529
Taiwan
Republic of China
e-mail: stcheng@stat.sinica.edu.tw